\documentclass[a4paper, 11pt]{article}

\usepackage[english]{babel}
\usepackage[latin2]{inputenc}
\usepackage[T1]{fontenc}

\usepackage{verbatim}
\usepackage[usenames]{xcolor}
\usepackage{soul}

\usepackage{amssymb}
\usepackage{amsmath}
\usepackage{amsthm}
	
\usepackage[overload]{empheq}

\usepackage[a4paper]{geometry}
\usepackage{caption}
\usepackage{graphicx}
\usepackage[font=normalsize]{caption}

\usepackage{setspace}
\usepackage{indentfirst}
\usepackage{lmodern}
\numberwithin{equation}{section}

\theoremstyle{plain}
\newtheorem{defi}{Definition}[section]
\newtheorem{theo}{Theorem}[section]
\newtheorem{prop}{Proposition}[section]

\newtheorem{lem}{Lemma}[section]
\newtheorem{cor}{Corollary}[section]
\newtheorem{rem}{Remark}[section]

\newcommand{\prove}{\noindent\textbf{Proof: }}
\newcommand{\quod}{\hfill$\square$}

\newcommand{\pt}{\partial_t}

\newcommand{\brck}[1]{{\left\vert\kern-0.25ex\left\vert\kern-0.25ex\left\vert #1 \right\vert\kern-0.25ex\right\vert\kern-0.25ex\right\vert}}

\usepackage{paralist}

\usepackage{enumerate}
\usepackage{enumitem}
\setenumerate{topsep=-2pt, itemsep=-1pt}

\usepackage[unicode, hidelinks, hyperfootnotes=false, pagebackref=false]{hyperref}

\usepackage{verbatim}

\begin{document}

\begin{center}

\textbf{\Large Wave-structure interaction for long wave models in the presence of a freely moving body on the bottom}\\
\vspace{0.3 cm}
\textbf{Krisztian Benyo}\footnote[2]{Institut de Math\'ematiques de Bordeaux, Universit\'e de Bordeaux, France. (\texttt{krisztian.benyo@math.u-bordeaux.fr})\vspace{0.1 cm}

The author was partially supported by the Conseil Régional Nouvelle Aquitaine, the ANR-13-BS01-0003-01 DYFICOLTI and the Fondation Del Duca.}\\

\end{center}

\vspace{1 cm}
\begin{center}\parbox{14cm}{
\noindent
\textsc{Abstract:} In this paper a particular fluid-solid interaction problem is addressed in which the solid object is lying at the bottom of a layer of fluid and moves under the forces created by waves traveling on the surface of this layer. More precisely, the water waves problem is considered in a fluid of fixed depth with a flat bottom topography and with an object lying on the bottom, allowed to move horizontally under the pressure forces created by the waves. After establishing the physical setting of the problem, namely the dynamics of the fluid and the mechanics of the solid motion, as well as analyzing the nature of the coupling, we examine in detail two particular shallow water asymptotic regimes: the case of the (nonlinear) Saint-Venant system, and the (weakly nonlinear) Boussinesq system. An existence and uniqueness theorem is proved for the coupled fluid-solid system in both cases. Using the particular structure of the coupling terms one is able to go beyond the standard scale for the existence time of solutions to the Boussinesq system with a moving bottom.}
\end{center}

\section*{Introduction}

The water waves problem, which consists in describing the motion of waves at the surface of an inviscid, incompressible, and irrotational fluid of constant density under the action of gravity, has attracted a lot of attention in the last decades. The local well-posedness theory is now well-understood following the works of Wu \cite{wuwellposedness2D, wuwellposedness3D} establishing the relevance of the Taylor sign condition. In the case of finite depth, which is of interest here, we refer for instance to \cite{lanneswp, iguchishallow, alazardgravitycauchy}; the case where the bottom is also allowed to depend on time has also been treated in \cite{alazardbottommove, iguchibottommove, melinandtsunami}. In this paper, we are interested in a particular configuration where the bottom depends on time, but instead of being in forced motion as in the above references, it evolves under the action of the hydrodynamic forces created by the surface waves. Finding its evolution is therefore a free boundary problem, which is coupled to the standard water waves problem, itself being a free boundary problem. The mathematical theory for such a configuration has not been considered yet; we refer however to \cite{lannesbonneton} for a related problem where the moving object is floating instead of lying on the bottom, as it is in the present paper.

Here, our goal is not to address the local well-posedness theory for this double free boundary problem, but to give some qualitative insight on its behavior by deriving and analyzing simpler asymptotic models. The focus is on a regime which is particularly interesting for applications, namely, the shallow water regime, where the typical horizontal scale of the flow is much larger than the depth at rest. For a fixed bottom, several models arise in this setting such as the Korteweg--de Vries (KdV) equation (justified in \cite{craiglangrangian,kanonishida,schneiderlongwavelimit}), the nonlinear shallow water equations (justified in \cite{ovsjannikov,kanonishida,alvarezlanneslong,iguchishallow}), the Boussinesq systems (justified in \cite{craiglangrangian,kanonishida,bonachen,davidboussinesq}) -- see also \cite{peregrineoriginal,papanicolaou,weikirbyboussinesq,chazelbottominfluence,craiglannesrough} for particular focus on topography effects --  the Green--Naghdi equations \cite{liapproximatif,alvarezlannesGN,iguchiGNtsunami}, etc. We refer to \cite{lannesbible} for more exhaustive references.

For a bottom with prescribed motion in general, the problem has already been considered, local well-posedness results (\cite{alazardbottommove}) and long time existence results (\cite{melinandtsunami}) have been proven recently. Numerical experiments and attempts to adapt existing and known shallow water models for a moving bottom regime have been present for a while in literature, however lacking rigorous justifications. After observing successively generated solitary waves due to a disturbance in the bottom topography advancing at critical speed (\cite{wusolitontrain}) they formally derived a set of generalized channel type Boussinesq systems (\cite{wutengwaterwaves}), their work was extended later on in a formal study on more general long wave regimes (\cite{chengeneralboussinesq}). Tsunami research has also proved to be an important motivating factor with the consideration of water waves type problems with a moving bottom (see for example \cite{guyennespectrum} or \cite{mitsotakistsunami} for an extensive numerical study). The mathematical justification of these models as approximations of the full problem was carried out not too long ago (\cite{iguchibottommove} for Saint-Venant type systems, or \cite{iguchiGNtsunami} for the precise Green--Naghdi system).

Here, we present a new class of problems where the bottom is still moving, but its movement is not prescribed, instead it is generated by the wave motion. A good approach to this is to place a freely moving object on the bottom of the fluid domain. The main physical motivation of this study lies in the recent development of submerged wave energy converters (submerged pressure differential devices, \cite{LehmannTheWC} and references therein) and oscillating wave surge converters (WaveRollers and Submerged plate devices, \cite{waveenergyconverters}), as well as reef-evolution and submarine landslide modeling problems. Bibliography in the more theoretical approach is rather lacking, existing studies are heavily oriented to physical experiments (most notably in \cite{abadibadou} where the authors investigate a submerged spring-block system and its numerical simulation through an adapted level set method, for further details, see for example \cite{cottetmaitrenumerics}), as well as numerical applications (\cite{zoazoazoo} for instance).

The structure of the article is as follows. In the first section the free surface fluid dynamics system and its possible reformulations in the water waves setting are presented. The equations governing the motion of an object lying on the bottom are established, they derive from Newton's second law and take into account the hydrodynamic force exerted by the fluid and a dynamic friction force. The characteristic scales of the variables of the system are also introduced in order to derive the nondimensionalized equivalents of the different equations and formulae, preparing for the study of the asymptotic models.

In Section 2, we detail the first order asymptotic regime with respect to the shallowness parameter $\mu$; the resulting approximation is the well-known (nonlinear) Saint-Venant equations, in the presence of a solid moving on the bottom of the fluid domain. A key step is to derive an asymptotic approximation of the hydrodynamic force exerted on the solid. Then we establish a local in time well-posedness result for the coupled system.

In the third section, we elaborate our study on a second order asymptotic regime with respect to the shallowness parameter $\mu$. This study concerns the so called long wave regime where the vertical size of the waves and of the solid are assumed to be small compared to the mean fluid height. The resulting approximation is the so called (weakly nonlinear) Boussinesq system. A local in time well-posedness is shown for this coupled system as well. The standard existence time for a Boussinesq system with a moving bottom is $\mathcal{O}(1)$ with respect to the nonlinearity parameter $\varepsilon$, due to the presence of a source term involving time derivatives of the topography, which can potentially become large (as remarked in \cite{melinandtsunami}). By a precise analysis of the wave-structure coupling one is able to extend the existence time to the $\mathcal{O}(\varepsilon^{-1/2})$ time scale. This time scale is therefore intermediate between the aforementioned $\mathcal{O}(1)$ scale, and the $\mathcal{O}(\varepsilon^{-1})$ scale that can be achieved for fixed bottoms (\cite{alvarezlanneslong, cosminboussinesq}).

\section{The fluid-solid coupled model}

In this section we present our model in general, that is, the equations characterizing the fluid dynamics as well as the equation describing the solid motion. The dimensionless equations are formulated at the end of the section to prepare the upcoming analysis for the shallow water asymptotic models.

\subsection{The dynamics of a fluid over a moving bottom}

As a basis for our model and our computations, we consider a fluid moving under the influence of gravity. The fluid domain $\Omega_t$ (depending on the time $t$) is delimited from below by a moving bottom and from above by a free surface. In our case the fluid is homogeneous with a constant density $\varrho$, moreover it is inviscid, incompressible, and irrotational.

To clarify the upcoming notations, the spatial coordinates take the form $(x,z)\in\mathbb{R}^d\times\mathbb{R}$ with $x$ denoting the horizontal component and $z$ the vertical one. Regarding differential operators, $x$ or $z$ as a subscript refers to the operator with respect to that particular variable, the absence of subscript for an operator depending on spatial variables means that it is to be taken for the whole space $(x,z)\in\mathbb{R}^{d+1}$. From a theoretical point of view arbitrary horizontal dimensions $d\in\mathbb{N}^+$ can be considered even though the physically relevant cases are $d=1$ and $2$ only.

In what follows, we denote by $\zeta(t,x)$ the free surface elevation function and $b(t,x)$ describes the bottom topography variation at a base depth of $H_0$. With this notation at our disposal the fluid domain is
\begin{equation*}
\Omega_t=\left\{(x,z)\in\mathbb{R}^d\times\mathbb{R}\,:\,-H_0+b(t,x)<z<\zeta(t,x)\right\},
\end{equation*}

Let us also introduce the height function $h(t,x)=H_0+\zeta(t,x)-b(t,x)$ that describes the total depth of the fluid at a given horizontal coordinate $x$ and at a given time $t$. 

In order to avoid special physical cases arising from the fluid domain $\Omega_t$ (such as islands or beaches), throughout our analysis we suppose the following (or similar) minimal water height condition
\begin{equation}\label{genminheightcond}
\exists h_{min}>0,\;\forall (t,x)\in[0,T)\times\mathbb{R}^d,\;h(t,x)\geqslant h_{min},
\end{equation}
we refer to \cite{poyferreshore} for an analysis of the water waves equation allowing vanishing depth, and to \cite{lannesmetivier} where the evolution of the shoreline is considered for the one dimensional nonlinear Saint-Venant and Serre--Green--Naghdi equations.

\subsubsection{The free surface Bernoulli equations}

To describe the fluid motion under the aforementioned physical assumptions, the free surface Euler equations could be considered, however for what follows the formulation involving a potential (the Bernoulli equations) is more adapted. Due to the fluid being incompressible and irrotational, one can describe its dynamics by utilizing the velocity potential $\Phi$, and with the knowledge of this potential one may recover the actual velocity field as the gradient.

The velocity potential is obtained as a solution of the following Laplace equation
\begin{equation}\label{eq_laplacegeneral}
\begin{cases}
\Delta\Phi = 0\quad\textnormal{ in }\Omega_t,\\
\Phi|_{z=\zeta}=\psi,\;\;\sqrt{1+|\nabla_xb|^2}\partial_{\mathbf{n}}\Phi|_{z=-H_0+b}=\pt b,
\end{cases}
\end{equation}
where $\psi$ is the velocity potential on the free surface (an unknown of the problem). Here we made use of the notation $\partial_{\mathbf{n}}$ signifying the upwards normal derivative (with $\mathbf{n}$ being the unit normal vector of the fluid domain pointing upward). Notice that the Neumann boundary condition on the bottom of the fluid domain corresponds to a kinematic (or no-penetration) boundary condition (that is, the fluid particles do not cross the bottom). Naturally the same condition applies to the free surface, meaning that
\begin{equation}\label{eq_bernoulliboundary}
\pt\zeta-\sqrt{1+|\nabla_x\zeta|^2}\partial_{\mathbf{n}}\Phi=0\quad \textnormal{ on }\{z=\zeta(t,x)\}.
\end{equation}

Additionally we also require that there is no surface tension along the free surface, so the pressure $P$ at the surface is given by the atmospheric pressure $P_{atm}$, hence
\begin{equation}\label{eq_eulerpressure}
P=P_{atm}\quad\textnormal{ on }\{z=\zeta(t,x)\}.
\end{equation}

By the momentum conservation of the fluid system we get that
\begin{equation}\label{eq_bernoullimoment}
\pt\Phi+\frac{1}{2}|\nabla\Phi|^2+gz = -\frac{1}{\varrho}(P-P_{atm})
\end{equation}
in the domain $\Omega_t$. Here $g$ in the equation denotes the gravitational acceleration, furthermore $\varrho$ denotes the density of the fluid (constant due to the homogeneity assumption).

So the free surface Bernoulli equations are the system of equations (\ref{eq_laplacegeneral})-(\ref{eq_bernoullimoment}).

Based on equation (\ref{eq_bernoullimoment}), we can recover the pressure in terms of the velocity potential:
\begin{equation}\label{form_pressioninit}
P=-\varrho\left(\pt\Phi+\frac{1}{2}|\nabla\Phi|^2+gz\right)+P_{atm}.
\end{equation}
This relation allows to compute the hydrodynamical force exerted on the solid by the fluid (derived from Newton's second law in Section \ref{sectsolid}).

\subsubsection{The Zakharov / Craig--Sulem framework}

We present another formulation of the equations (also referred to as the water waves problem). This formulation is attributed to Zakharov in his studies regarding gravity waves \cite{zakharovformula} and is based on the fact that the variables $\zeta$ and $\psi=\Phi|_{z=\zeta}$ fully determine the flow. More precisely, the water waves problem reduces to a set of two evolution equations in $\zeta$ and $\psi$,
\begin{equation}\label{eq_zakharov}
\begin{cases}
\pt\zeta-\sqrt{1+|\nabla_x\zeta|^2}\partial_\mathbf{n}\Phi|_{z=\zeta} = 0,\\
\pt\psi+g\zeta+\dfrac{1}{2}|\nabla_x\psi|^2-\dfrac{(\sqrt{1+|\nabla_x\zeta|^2}\partial_\mathbf{n}\Phi|_{z=\zeta}+\nabla_x\zeta\cdot\nabla_x\psi)^2}{2(1+|\nabla_x\zeta|^2)}=0,
\end{cases}
\end{equation}
where $\Phi$ solves the boundary value problem (\ref{eq_laplacegeneral}).

In more general terms, one can introduce a natural decomposition of $\Phi$ into a ``fixed bottom'' and a ``moving bottom'' component which could be used to define the so-called Dirichlet-Neumann and Neumann-Neumann operators associated to the Laplace problem (\ref{eq_laplacegeneral}) (assuming sufficient regularity for the limiting functions), but we will not pursue further this path, for more details we refer to the works of Craig and Sulem \cite{craigsulem2, craigsulem1}. For a more specific analysis of the moving bottom case we refer to the article of Alazard, Burq, and Zuily \cite{alazardbottommove} for the local well-posedness theory or to \cite{iguchibottommove} for specific studies motivated by earthquake generated tsunami research. For a comprehensive and detailed analysis as well as the well-posedness of the water waves problem in the general setting, we refer to \cite{lannesbible} and references therein.

Since our study focuses on shallow water regimes, it is convenient to bypass the aforementioned technicalities by introducing the following variable:
\begin{defi} The vertically averaged horizontal component of the velocity is given by
\begin{equation}\label{def_averagedvelocity}
\overline{V} = \frac{1}{h}\int_{-H_0+b}^{\zeta}\nabla_x\Phi(\cdot,z)\,dz,
\end{equation}
where $\Phi$ solves (\ref{eq_laplacegeneral}).
\end{defi}
The interest of this new variable $\overline{V}$ is that a closed formulation of the water waves problem in terms of $\zeta$ and $\overline{V}$ (instead of $\zeta$ and $\psi$) can be obtained, see for example \cite{lannesbonneton}. For our case, it is sufficient to observe that
\begin{prop}\label{prop_nonasymptoticrewrite}
If $\Phi$ solves (\ref{eq_laplacegeneral}) and $\overline{V}$ is defined as in (\ref{def_averagedvelocity}), then
\begin{equation}\label{dirichletneumannrel}
\sqrt{1+|\nabla_x\zeta|^2}\partial_{\mathbf{n}}\Phi|_{z=\zeta}=\partial_tb-\nabla\cdot(h\overline{V}),
\end{equation}
assuming sufficient regularity on the data concerning $\zeta$, $\psi$, and $b$ as well as the minimal water depth condition (\ref{genminheightcond}).
\end{prop}
\begin{rem}
Let $\zeta,b\in W^{1,\infty}(\mathbb{R}^d)$ such that they satisfy the minimal water depth condition (\ref{genminheightcond}). Moreover, let $\psi\in\dot{H}^{3/2}(\mathbb{R}^d)=\{f\in L^2_{loc}\,:\,\nabla_x f\in H^{1/2}(\mathbb{R}^d)\}$. Then the Laplace equation (\ref{eq_laplacegeneral}) can be solved with $\Phi\in\dot{H}^2(\Omega)=\{f\in L^2_{loc}\,:\,\nabla_x f\in H^1(\Omega)\}$ and relation (\ref{dirichletneumannrel}) holds true, where
\begin{equation}
\Omega = \{(X,z)\in\mathbb{R}^d\times\mathbb{R},\;-1+b(X)<z<\zeta(X)\}
\end{equation}
is a known fluid domain. For more details, we refer to Chapter 2 of \cite{lannesbible}.
\end{rem}
With this, the water waves problem with a moving bottom takes the following form
\begin{equation}\label{eq_zakharovmovbott}
\begin{cases}
\pt\zeta+\nabla\cdot(h\overline{V}) = \partial_tb,\\
\pt\psi+g\zeta+\dfrac{1}{2}|\nabla_x\psi|^2-\dfrac{(-\nabla\cdot(h\overline{V})+\partial_tb+\nabla_x\zeta\cdot\nabla_x\psi)^2}{2(1+|\nabla_x\zeta|^2)}=0.
\end{cases}
\end{equation}
This system seemingly depends on three variables, namely $\zeta$, $\overline{V}$ and $\psi$, but in fact the Laplace equation provides a connection between the latter two. Exploiting this connection to express (asymptotically) one variable with the other gives rise to various well-known asymptotic equations under the shallow water assumption. In Section \ref{sec_nodim} detailing the nondimensionalization of the system we shall provide the necessary tools as well as some references concerning this asymptotic expansion.

\subsection{A freely moving object on a flat bottom}\label{sectsolid}

The aim of this paper is to understand a particular case in which the bottom of the domain contains a freely moving object, the movement of which is determined by the gravity driven fluid motion. We will work with a flat bottom in the presence of a freely moving solid object on it (see Figure \ref{fig_domainwithsolid}).

\begin{center}
\captionsetup{type=figure}
\includegraphics[width=0.8\linewidth]{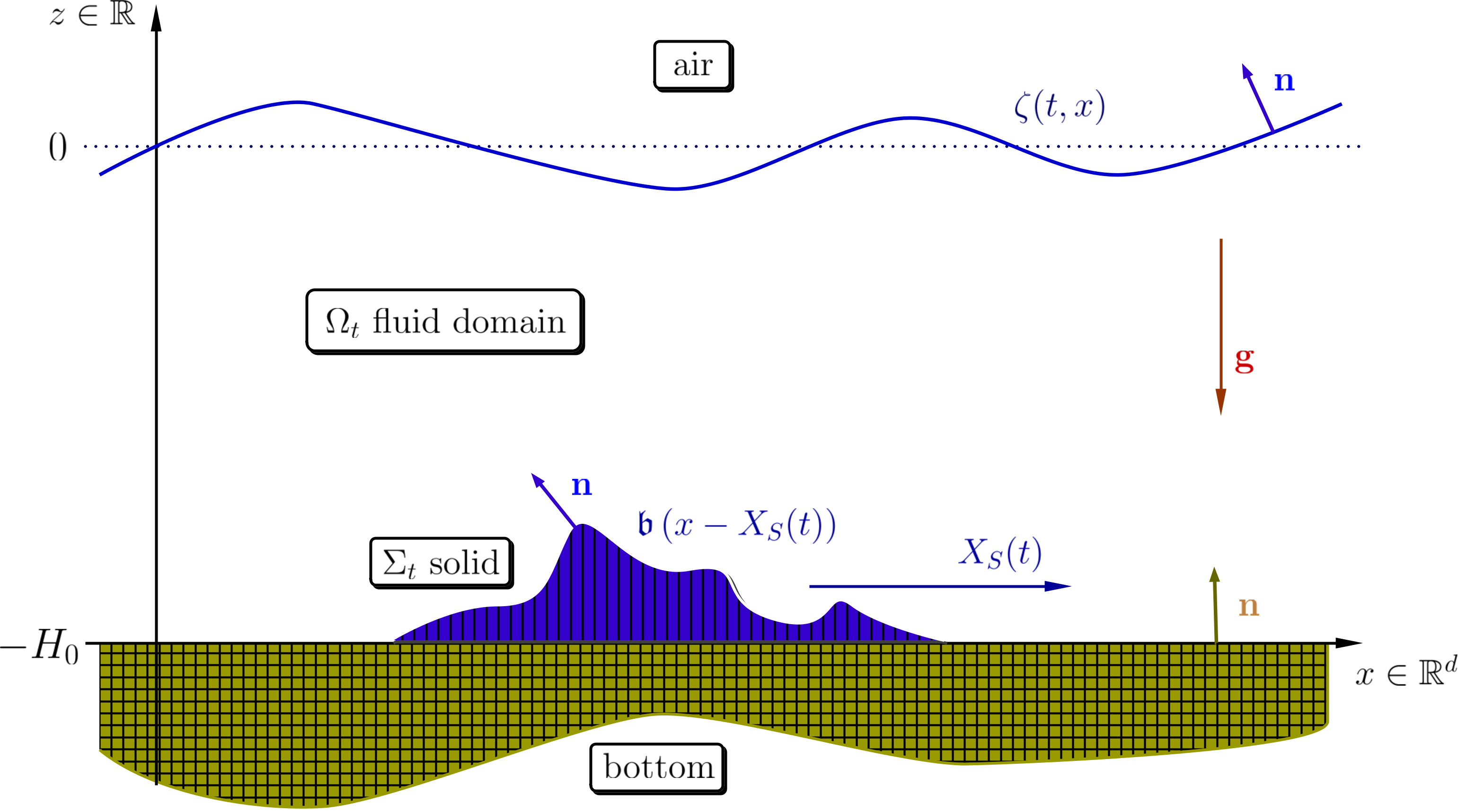}
\captionof{figure}{The setting of the water waves problem in the presence of a solid on the bottom}\label{fig_domainwithsolid}
\end{center}

For the solid we suppose it to be rigid and homogeneous with a given mass $M$. The surface of the object can be characterized by two components: the part of the surface in direct contact with the fluid, denoted by $\Sigma_t$ and the rest, that is the part in direct contact with the flat bottom, denoted by $I(t)$. For convenience reasons we shall suppose that $\Sigma_t$ is a graph of a $\mathcal{C}^{\infty}$ function with compact support $I(t)$ for any instance of $t$.

The solid moves horizontally in its entirety, we denote by $X_S(t)$ the displacement vector, and $v_S(t)$ the velocity (with $\dot{X}_S=v_S$). We make the additional hypothesis that the object is neither overturning, nor rotating so its movement is completely described by its displacement vector, which will be restrained to horizontal movement only. In particular, this means that the object is not allowed to start floating, the domain $I(t)$ has a constant (nonzero) area.

Under these assumptions a simplified characterization of the function describing the bottom variation is possible:
\begin{equation}\label{form_bottomfunction}
b(t,x)=\mathfrak{b}\left(x-X_S(t)\right),
\end{equation}
where $\mathfrak{b}$ corresponds to the initial state of the solid at $t=0$ (so that we have $X_S(0)=0$).

Taking into account all the external forces acting on the object, Newton's second law provides us with the correct equation for the movement of the solid. The total force acting on the solid is
\begin{equation*}\begin{aligned}
\mathbf{F}_{total} &= \mathbf{F}_{gravity} + \mathbf{F}_{solid-bottom\; interaction} +  \mathbf{F}_{solid-fluid\; interaction}\\
&= M\mathbf{g} + \left[\mathbf{F}_{normal} +\mathbf{F}_{friction}\right]+ \mathbf{F}_{pressure}.
\end{aligned}\end{equation*}
Here we made use of the fact that the force emerging from the contact of the solid with the bottom may be decomposed in two components: the normal force, perpendicular to the surface of the bottom, expressing the fact that the bottom is supporting the solid, and the (kinetic or dynamic) friction force, the tangential component, hindering the sliding of the solid. By making use of the three empirical laws of friction \cite{berthelotmechanic}, most notably the third law often attributed to Coulomb regarding the existence of a coefficient $c_{fric}>0$ of kinetic friction (describing the material properties of the contact medium), we may reformulate the tangential contact force as follows
\begin{equation}\label{form_frictionlaw}
\mathbf{F}_{friction}=\mathbf{F}_{sliding\; friction}=-c_{fric}|\mathbf{F}_{normal}|\frac{v_S(t)}{|v_S(t)|+\sqrt{gH_0}\overline{\delta}},
\end{equation}
where $\overline{\delta}\ll 1$ is a purely mathematical dimensionless parameter serving as a regularizing term in order to avoid a singularity in the equation when the solid stops, that is when $v_S(t)$ is equal to $0$. Normally, when the solid comes to a halt, the kinetic friction detailed just before turns into static friction, a tangential force component preventing the solid from restarting its movement. The static friction has its own coefficient, which is usually greater than $c_{fric}$, and its direction is determined by the horizontal force component rather than the velocity. 

\begin{rem} 
The coefficient of friction $c_{fric}$ is a dimensionless scalar constant, it describes a ratio proportional to the hindering effect generated by the parallel motion of two surfaces. It is in fact a property of the system, in reality it not only depends on the material of the two surfaces but their geometry (surface microstructure), temperature, atmospheric conditions, velocity of the motion, etc. and as such it is impossible to accurately determine it. To give the reader an idea, an almost frictionless sliding (for example objects on ice, lubricated materials) corresponds to a coefficient of $10^{-2}\sim 10^{-3}$, while a frictional sliding (for example rubber on paper) has a coefficient of order $1$.
\end{rem}

In order to prevent the complications that would arise by implementing the physically more relevant threshold for $v_S(t)=0$ and the associated jump in friction force, we simplify the system by regularizing the friction force, thus neglecting static effects. A more specific modeling and analysis of the transition between static and dynamic friction will be addressed in future works. 

Treating the horizontal and vertical component of $\mathbf{F}_{total}=(\mathbf{F}_{total}^h,\mathbf{F}_{total}^v)^\top$ separately and using the fact that the solid is constrained to horizontal motion, we have that the vertical components are in equilibrium, thus
\begin{equation}\label{eq_verticalforceequi}
0 = -Mg + \mathbf{F}_{normal}+\mathbf{F}_{pressure}^v,
\end{equation}
and we obtain that the horizontal movement of the solid is given by
\begin{equation}\label{eq_horizontalforceequi}
M\ddot{X}_S(t)=\mathbf{F}_{sliding\; friction}+ \mathbf{F}_{pressure}^h.
\end{equation}

Finally, by making use of the fact that 
\begin{equation*}
\mathbf{F}_{pressure}=\int_{\Sigma_t}P\mathbf{n}_{solid}\,d\Sigma=\int_{I(t)}P|_{z=-H_0+b(t,x)}\begin{pmatrix}\nabla_xb\\-1\end{pmatrix}\,dx,
\end{equation*}
due to the fact that the inwards normal vector for the surface of the solid $\mathbf{n}_{solid}=-\mathbf{n}$ can be easily expressed by the bottom variation $b(t,x)$, since
\begin{equation*}
\mathbf{n}_{solid}=\frac{1}{\sqrt{1+|\nabla_xb|^2}}\begin{pmatrix} \nabla_xb\\ -1\end{pmatrix}.
\end{equation*}

Therefore we obtain from (\ref{eq_verticalforceequi}) that
\begin{equation}
F_{normal}=Mg+\int_{I(t)}P|_{z=-H_0+b(t,x)}\,dx,
\end{equation}
now as a scalar quantity since the vertical direction is one dimensional. Therefore, by (\ref{form_frictionlaw}), (\ref{eq_horizontalforceequi}) writes as
\begin{equation}
M\ddot{X}_S(t)=-c_{fric}F_{normal}\frac{\dot{X}_S(t)}{\left|\dot{X}_S(t)\right|+\sqrt{gH_0}\overline{\delta}}+\int_{I(t)}P|_{z=-H_0+b(t,x)}\nabla_xb\,dx.
\end{equation}

So we have that Newton's equation characterizing the motion of the solid takes the following form
\begin{equation}\label{eq_newton}
M\ddot{X}_S(t)=-c_{fric}\left(Mg+\int_{I(t)}P|_{z=-H_0+b(t,x)}\,dx\right)\frac{\dot{X}_S(t)}{\left|\dot{X}_S(t)\right|+\sqrt{gH_0}\overline{\delta}}+\int_{I(t)}P|_{z=-H_0+b(t,x)}\nabla_xb\,dx.
\end{equation}

A key step in our study is to handle the force term exerted by the fluid, which requires the computation of the integral of the pressure on the bottom over the solid domain. For this we will establish an appropriate formula for the pressure to be used in the integral.

In both the case of the freely moving bottom (due to the moving object) and the free surface, the kinematic no-penetration condition still applies, most notably we still have that
\begin{equation*}
\partial_tb - \sqrt{1+|\nabla_xb|^2}\mathbf{U}\cdot \mathbf{n} = 0 \quad\textnormal{ for } \{z=-H_0+b(t,x)\},
\end{equation*}
or equivalently, on the part of the surface of the solid in contact with the fluid ($\Sigma_t$), the normal component of the fluid velocity field coincides with the normal component of the velocity of the solid, that is
\begin{equation}
\mathbf{U}\cdot \mathbf{n}_{solid} = v_S\cdot \mathbf{n}^h_{solid} \quad\textnormal{ for } \{z=-H_0+b(t,x)\}.
\end{equation}

To sum up, the water waves problem in the presence of a solid on the bottom is given by equations (\ref{eq_zakharovmovbott}) and (\ref{eq_laplacegeneral}), where in the Neumann boundary condition, the bottom function $b$ and its time derivative are given by (\ref{form_bottomfunction}), with $X_S$ arising from (\ref{eq_newton}) and the pressure $P$ derived from (\ref{form_pressioninit}).

\subsection{Dimensionless form of the equations}\label{sec_nodim}

The main part of the analysis consists of establishing and analyzing the wave-structure interaction system for shallow water regimes, for that we need first of all the correct parameters involving the characteristic orders of magnitude of our variables as well as the dimensionless equations obtained with the help of these quantities.

\subsubsection{The different scales of the problem}

 First of all we present the proper dimensionless parameters relevant to the system. For that we need to introduce the various characteristic scales of the problem: as already mentioned before, the base water depth is $H_0$. The characteristic horizontal scale of the wave motion (both for longitudinal and transversal directions) is $L$, the order of the free surface amplitude is $a_{surf}$, and the characteristic height of the solid (order of the bottom topography variation in general) is $a_{bott}$.

Using these quantities, we can introduce several dimensionless parameters:
\begin{equation*}
\mu=\frac{H_0^2}{L^2},\quad \varepsilon=\frac{a_{surf}}{H_0},\quad\textnormal{ and }\beta=\frac{a_{bott}}{H_0},
\end{equation*}
where $\mu$ is called the shallowness parameter, $\varepsilon$ stands for the nonlinearity (or amplitude) parameter, and $\beta$ is the bottom topography parameter.

Our goal in this paper is to examine asymptotic models when $\mu$ is small (shallow water regime), and under various assumptions on the characteristic size of $\varepsilon$ and $\beta$.

With these parameters in our hand, we may remark that the natural scaling for the horizontal spatial variable $x$ is $L$, and for its vertical counterpart $z$ it is $H_0$. Moreover the natural order of magnitude for the function characterizing the free surface $\zeta$ is $a_{surf}$, and for the bottom $b$ it is $a_{bott}$. Thus the nondimensionalized form for the water depth is
\begin{equation*}
h=1+\varepsilon\zeta-\beta b.
\end{equation*} 

Furthermore, one can establish the correct scale of the velocity potential through linear wave analysis, which gives rise to
\begin{equation*}
\Phi_0=\frac{a_{surf}}{H_0}L\sqrt{gH_0}.
\end{equation*}
As for the pressure, we choose the typical order of the hydrostatic pressure, that is $P_0=\varrho gH_0$. For the time parameter, from linear wave theory one can deduce the scaling as
\begin{equation*}
t_0=\frac{L}{\sqrt{gH_0}}.
\end{equation*}

Finally, for the parameters concerning the solid, we impose that the characteristic horizontal dimension of the solid is comparable to $L$ (which was already implicitly assumed). It would be relevant to consider solids with a smaller size, but this raises important difficulties. Even in the case of a fixed bottom there is no fully justified model yet in general (see for example \cite{craiglannesrough}).

Following this, by taking into account the volume integral of the density, the natural scaling for the mass of the solid is given by
\begin{equation*}
M=M_0\tilde{M}=\varrho L^da_{bott}\tilde{M}.
\end{equation*}

Thus the proper nondimensionalized parameters are obtained by
\begin{equation*}
x'=\frac{x}{L},\quad z'=\frac{z}{H_0},\quad\zeta'=\frac{\zeta}{a_{surf}},\quad\Phi'=\frac{\Phi}{\Phi_0},\quad t'=\frac{t}{t_0},\;\;\textnormal{etc.}
\end{equation*}
For the sake of clarity we shall omit the primes on the variables from here on.

Our main interest will be to express the equations principally with the different orders of magnitude of $\mu$ (the shallowness parameter) to pass on to the different asymptotic regimes. Given the particular structure of the asymptotic regimes we are going to examine we shall make an a priori hypothesis concerning certain parameters.
\begin{rem}\label{rem_blessthane}
Since all the regimes handled in this article involve the hypothesis that $\varepsilon$ and $\beta$ are of the same order of magnitude we assume, for the sake of simplicity, that $\beta=\varepsilon$.
\end{rem}

An additional precision shall be made concerning the quantities involving the bottom. The explicit form of the nondimensionalized form for the water depth is
\begin{equation}\label{eq_nondimheight}
h(t,x)=1+\varepsilon(\zeta(t,x)-b(t,x))
\end{equation} 
with
\begin{equation}\label{form_ndbottomvari}
b(t,x)=\mathfrak{b}\left(x-X_S(t)\right).
\end{equation}

\subsubsection{Nondimensionalized equations}

Using the previous section and in particular taking $\varepsilon=\beta$ as in Remark \ref{rem_blessthane}, one easily derives the dimensionless version of (\ref{eq_zakharovmovbott}), namely
\begin{equation}\label{eq_ndzakharovmovingbottom}
\begin{cases}
\pt\zeta+\nabla\cdot(h\overline{V}) = \partial_tb,\\
\pt\psi+\zeta+\dfrac{\varepsilon}{2}|\nabla_x\psi|^2-\varepsilon\mu\dfrac{(-\nabla\cdot(h\overline{V})+\partial_tb+\nabla_x(\varepsilon\zeta)\cdot\nabla_x\psi)^2}{2(1+\varepsilon^2\mu|\nabla_x\zeta|^2)}=0,
\end{cases}
\end{equation}
where $\overline{V}$ is now defined as
\begin{equation}
\overline{V} = \frac{1}{h}\int_{-1+\varepsilon b}^{\varepsilon\zeta}\nabla_x\Phi(\cdot,z)\,dz,
\end{equation}
with $h=1+\varepsilon\zeta-\varepsilon b$, furthermore $\Phi$ solves
\begin{equation}\label{eq_nondimlaplace}
\begin{cases}
\Delta^{\mu}\Phi = \mu\Delta_x\Phi + \partial_z^2\Phi = 0,\quad \textnormal{on }-1+\varepsilon b\leqslant z\leqslant\varepsilon\zeta,\\
\Phi|_{z=\varepsilon\zeta}=\psi,\quad \left(\partial_z\Phi-\mu\nabla_x(\varepsilon b)\cdot\nabla_x\Phi\right)|_{z=-1+\varepsilon b}=\mu\partial_tb,
\end{cases}
\end{equation}
the nondimensionalized equivalent of the Laplace problem (\ref{eq_laplacegeneral}).

It is also necessary to nondimensionalize the formula describing the pressure (\ref{form_pressioninit}), thus
\begin{equation}\label{form_pressurewithphi}
P=\frac{P_{atm}}{\varrho gH_0}-z-\varepsilon\partial_t\Phi-\frac{\varepsilon^2}{2}|\nabla_x\Phi|^2-\frac{\varepsilon^2}{2\mu}|\partial_z\Phi|^2.
\end{equation}
Here we had to separate the horizontal and the vertical part of the gradient due to the different scaling parameters for the different directions.

We remark that the normal derivative is given by
\begin{equation*}
\mathbf{n}_{solid}=\frac{1}{\sqrt{1+\varepsilon^2\mu|\nabla_xb|^2}}\begin{pmatrix} \sqrt{\mu}\varepsilon\nabla_xb\\ -1\end{pmatrix}.
\end{equation*}
Thus we may reformulate Newton's equation (\ref{eq_newton}) in the following way
\begin{equation}\label{eq_ndnewton}
\ddot{X}_S(t)
=-\frac{c_{fric}}{\sqrt{\mu}}\left(1+\frac{1}{\varepsilon\tilde{M}}\int_{I(t)}P|_{z=-1+\varepsilon b}\,dx\right)\frac{\dot{X}_S(t)}{\left|\dot{X}_S(t)\right|+\overline{\delta}}+\frac{1}{\tilde{M}}\int_{\mathbb{R}^d}P_{z=-1+\varepsilon b}\nabla_{x}b\,dx,
\end{equation}
taking into consideration the characteristic scales of the variables.

\section{The $\mathcal{O}(\mu)$ asymptotic regime: The nonlinear Saint-Venant equations}\label{sec_saintvenant} 

We shall now start our analysis for shallow water regimes, that is an asymptotic analysis with respect to the shallowness parameter $\mu$ for the nondimensionalized water waves problem (\ref{eq_ndzakharovmovingbottom}) coupled with Newton's equation (\ref{eq_ndnewton}) for the solid. With our notations, this means that we would like to consider systems that are valid for $\mu\ll 1$. 

In this section we treat the general first order approximate system, more specifically a model with $\mathcal{O}(\mu)$ approximation that allows large wave amplitudes and large bottom variations ($\varepsilon=\mathcal{O}(1)$). So, the asymptotic regime writes as follows
\begin{equation}\tag{SV}\label{hyp_saintvenantregime}
0\leqslant\mu\leqslant\mu_{max}\ll 1,\quad\varepsilon=1.
\end{equation}

\subsection{The fluid equations in the asymptotic regime}

As mentioned before, the important step in deducing asymptotic models relies on how we establish the connection between the variables $\overline{V}$ and $\psi$. More precisely, it is possible to construct an asymptotic expansion of $\overline{V}$ with respect to $\mu$ (depending on $\zeta$, $b$ and $\psi$). For details, we refer to Chapter $3$ of \cite{lannesbible}. One can equally obtain an asymptotic expansion of $\Phi$ with respect to $\mu$, depending on the aforementioned variables.
Quite obviously the equation $\Delta^\mu\Phi = 0$ in (\ref{eq_nondimlaplace}) reduces to $\partial_z^2\Phi = 0$ at leading order in $z$; since the Neumann boundary condition in (\ref{eq_nondimlaplace}) is $\mathcal{O}(\mu)$, it follows that $\Phi$ does not depend on $z$ at leading order, and therefore
\begin{equation*}
\overline{V}=\nabla_x\psi+\mathcal{O}(\mu),
\end{equation*}
see Proposition 3.37. in \cite{lannesbible} for a rigorous proof.

So the system (\ref{eq_ndzakharovmovingbottom}) for the $(\zeta,\,\overline{V})$ variables simplifies as follows
\begin{equation}\label{eq_nonlinearSV}
\begin{cases}
\partial_t\zeta+\nabla_x\cdot(h\overline{V}) = \partial_t b,\\
\partial_t\overline{V}+\nabla_x\zeta+(\overline{V}\cdot\nabla_x)\overline{V} = 0,
\end{cases}
\end{equation}
where we considered the gradient of the second equation in (\ref{eq_ndzakharovmovingbottom}), and then neglected terms of order $\mathcal{O}(\mu)$. This system is known as the (nonlinear) Saint-Venant or nonlinear shallow water equations.

\subsection{Formal derivation of a first order asymptotic equation for the solid motion}\label{sec_derivationnewton1}

Our strategy is as follows: we establish an asymptotic formula of order $\mathcal{O}(\mu)$ for the pressure $P$ based on (\ref{form_pressurewithphi}). With this at our disposal, we rewrite Newton's equation (\ref{eq_ndnewton}) at order approximately $\mu$ describing the displacement of the solid.

For an $\mathcal{O}(\mu)$ approximation, we shall start with the corresponding development for the velocity potential, that is
\begin{equation}\label{form_asymptoticphi1}
\Phi = \psi + \mathcal{O}(\mu),
\end{equation}
where $\psi=\Phi|_{z=\varepsilon\zeta}$ as before, the restriction of the velocity potential on the free surface. Knowing this we recover the following for the time derivative of $\psi$ (based on the second equation of the water waves problem (\ref{eq_ndzakharovmovingbottom}))
\begin{equation*}
\pt\psi = -\zeta - \frac{1}{2}|\nabla_x\psi|^2 + \mathcal{O}(\mu). 
\end{equation*}

So by substituting the first order asymptotic expansion of the velocity potential described in (\ref{form_asymptoticphi1}) into the general nondimensionalized formula of the pressure (\ref{form_pressurewithphi}) the corresponding $\mathcal{O}(\mu)$ approximation for the pressure takes the form
\begin{equation*}
P=\frac{P_{atm}}{\varrho gH_0}+(\zeta-z)+\mathcal{O}(\mu),
\end{equation*}
using the fact that $\psi$ does not depend on the variable $z$.

So in particular, at the bottom, we find that the pressure is given by the hydrostatic formula
\begin{equation}\label{form_pressureassymptotic}
P|_{z=-1+b}=\frac{P_{atm}}{\varrho gH_0}+h+\mathcal{O}(\mu).
\end{equation}

Thus for Newton's equation (\ref{eq_ndnewton}),
\begin{equation*}
\begin{aligned}
\ddot{X}_S &= -\frac{c_{fric}}{\sqrt{\mu}}\left(1+\frac{1}{\tilde{M}}\int_{I(t)}\left(\frac{P_{atm}}{\varrho gH_0}+h\right)\,dx\right)\frac{\dot{X}_S}{\left|\dot{X}_S\right|+\overline{\delta}}\\
&\qquad+\frac{P_{atm}}{\varrho gH_0\tilde{M}}\int_{\mathbb{R}^d}\nabla_xb\,dx+\frac{1}{\tilde{M}}\int_{\mathbb{R}^d}h\nabla_xb\,dx+\mathcal{O}\left(\frac{c_{fric}}{\tilde{M}}\sqrt{\mu}\right).
\end{aligned}
\end{equation*}

Using the fact that $b$ is of compact support, the integral of its (and $b^2$'s) gradient on the whole horizontal space is $0$, and the equation simplifies to
\begin{equation*}
\begin{aligned}
\ddot{X}_S &= -\frac{c_{fric}}{\sqrt{\mu}}\left(1+\frac{|\operatorname{supp}(\mathfrak{b})|}{\tilde{M}}\left(\frac{P_{atm}}{\varrho gH_0}+1\right)-\frac{|\operatorname{Volume}_{Solid}|}{\tilde{M}}+\frac{1}{\tilde{M}}\int_{I(t)}\zeta\,dx\right)\frac{\dot{X}_S}{\left|\dot{X}_S\right|+\overline{\delta}} \\
&\qquad+\frac{1}{\tilde{M}}\int_{\mathbb{R}^d}\zeta\nabla_xb\,dx+\mathcal{O}\left(\frac{c_{fric}}{\tilde{M}}\sqrt{\mu}\right).
\end{aligned}
\end{equation*}

Notice the presence of the friction term (the first term on the right hand side). Even though it is of order $\mu^{-1/2}$, it will not pose a problem when controlling the solid velocity, as we are going to see in Lemma \ref{lemma_velocityestimate}. later on (since it acts as a damping force). 

Recalling that $b$ is given by (\ref{form_ndbottomvari}) the corresponding approximative equation characterizing the motion of the body is
\begin{equation}\label{eq_newtonorder1}
\ddot{X}_S=-\frac{c_{fric}}{\sqrt{\mu}}\left(c_{solid}+\frac{1}{\tilde{M}}\int_{\operatorname{supp}(\mathfrak{b})+X_S}\zeta\,dx\right)\frac{\dot{X}_S}{\left|\dot{X}_S\right|+\overline{\delta}}+\frac{1}{\tilde{M}}\int_{\mathbb{R}^d}\zeta\nabla_x\mathfrak{b}(x-X_S)\,dx,
\end{equation}
where we made use of the following abbreviation:
\begin{equation}\label{form_solidcoefficient}
c_{solid}=1+\frac{|\operatorname{supp}(\mathfrak{b})|}{\tilde{M}}\left(\frac{P_{atm}}{\varrho gH_0}+1\right)-\frac{|\operatorname{Volume}_{Solid}|}{\tilde{M}}.
\end{equation}

Therefore, we have the following concerning the consistency of the solid equation:
\begin{prop}
Let $s_0\geqslant 0$, and let us assume that $\zeta\in\mathcal{C}([0,T];H^{s_0+4}(\mathbb{R}^d))$ and that $\mathfrak{b}\in H^{s_0+4}(\mathbb{R}^d)$ compactly supported. Furthermore let us suppose that $\nabla_x\psi\in \mathcal{C}([0,T];H^{s_0+4}(\mathbb{R}^d))$. The solid equation (\ref{eq_ndnewton}) is consistent at order $\mathcal{O}(\sqrt{\mu})$ with the model (\ref{eq_newtonorder1}) on $[0,T]$ with $T>0$.
\end{prop}
\prove By the regularity assumptions (Lemma 3.42. of \cite{lannesbible}), we can write that $\Phi = \psi + \mu R_1$ with
\begin{equation*}
\begin{aligned}
\brck{R_1}_{T,H^{s_0}}&\leqslant C(\brck{\zeta}_{T,H^{s_0+2}}, \|\mathfrak{b}\|_{H^{s_0+2}})\brck{\nabla_x\psi}_{T,H^{s_0+2}}\\
\brck{\partial_tR_1}_{T,H^{s_0}}&\leqslant C(\brck{\zeta}_{T,H^{s_0+4}}, \|\mathfrak{b}\|_{H^{s_0+4}},\brck{\nabla_x\psi}_{T,H^{s_0+4}}),
\end{aligned}
\end{equation*}
here the second estimate is due to Lemma 5.4. of \cite{lannesbible}. This means that, following the same computations as before, we have that
\begin{equation*}
\begin{aligned}
P|_{z=-1+b}&=\frac{P_{atm}}{\varrho gH_0}+h+\mu R_{P,1}\\
\brck{R_{P,1}}_{T,H^{s_0}}&\leqslant C(\brck{\zeta}_{T,H^{s_0+4}}, \|\mathfrak{b}\|_{H^{s_0+4}},\brck{\nabla_x\psi}_{T,H^{s_0+4}}).
\end{aligned}
\end{equation*}
Here the $\brck{ }_{T,\mathcal{X}}$ notation was adopted based on Definition \ref{def_bracket}.

Hence, in the equation for the solid motion (\ref{eq_ndnewton}), we recover the approximate equation (\ref{eq_newtonorder1}) with the additional error terms
\begin{equation*}
-\sqrt{\mu}\frac{c_{fric}}{\tilde{M}}\frac{\dot{X}_S}{\left|\dot{X}_S\right|+\overline{\delta}}\int_{I(t)}R_{P,1}\,dx+\mu\frac{1}{\tilde{M}}\int_{I(t)}R_{P,1}\nabla_x\mathfrak{b}(x-X_S)\,dx,
\end{equation*}
that can be estimated as an $\mathcal{O}(\sqrt{\mu})$ total error term, that is, it is less than
\begin{equation*}
\sqrt{\mu} C(\tilde{M}^{-1}, \brck{\zeta}_{T,H^{s_0+4}}, \|\mathfrak{b}\|_{H^{s_0+4}},\brck{\nabla_x\psi}_{T,H^{s_0+4}}).
\end{equation*}\quod

\subsection{The wave-structure interaction problem at first order}

With (\ref{eq_newtonorder1}) in our hand, we have all three equations for our coupled system. Indeed, notice that for the first equation in the nonlinear Saint-Venant system (\ref{eq_nonlinearSV}), the right hand side depends on $X_S$, since $b(t,x)$ depends on it. Hence, by the chain rule the right hand side is
\begin{equation*}
\pt b(t,x)=-\nabla_x\mathfrak{b}(x-X_S(t))\cdot\dot{X}_S(t).
\end{equation*}
Our remark concerning the friction term present in the acceleration equation (\ref{eq_newtonorder1}) becomes even more pertinent now, since we can observe a direct influence of the solid velocity in the first equation of the fluid system (\ref{eq_nonlinearSV}). This implies that a careful attention has to be paid on the velocity estimate for the solid.

\textit{To sum it up, the free surface equations with a solid moving at the bottom in the case of the nonlinear Saint-Venant approximation take the following form}
\begin{subequations}\label{eq_systemorder1}
\begin{align}[left = {\,\empheqlbrace}]
\begin{split}
&\partial_t\zeta+\nabla_x\cdot(h\overline{V}) = \nabla_x\mathfrak{b}(x- X_S)\cdot\dot{X}_S,\\
&\partial_t\overline{V}+\nabla_x\zeta+(\overline{V}\cdot\nabla_x)\overline{V} = 0,
\end{split} \label{eq_systemorder1o1} \\ 
&\ddot{X}_S = -\frac{c_{fric}}{\sqrt{\mu}}\left(c_{solid}+\frac{1}{\tilde{M}}\int_{I(t)}\zeta\,dx\right)\frac{\dot{X}_S}{\left|\dot{X}_S\right|+\overline{\delta}}+\frac{1}{\tilde{M}}\int_{\mathbb{R}^d}\zeta\nabla_x\mathfrak{b}(x-X_S)\,dx. \label{eq_systemorder1o2}
\end{align}
\end{subequations}

In what follows, we proceed to the mathematical analysis of this system. We shall establish a local in time existence result for the coupled equations.

\subsection{Local in time existence of the solution}

The main result on the local well-posedness of the wave-structure interaction problem (\ref{eq_systemorder1}) is the following:

\begin{theo}
Suppose that $\varepsilon=1$, and that $\mu$ is sufficiently small so that we are in the shallow water regime (\ref{hyp_saintvenantregime}). Let us suppose that for the initial value $\zeta_{in}$ and $\mathfrak{b}$ the lower bound condition (\ref{genminheightcond}) is satisfied. If the initial values $\zeta_{in}$ and $\overline{V}_{in}$ are in $H^s(\mathbb{R}^d)$ with $s\in\mathbb{R}$, $s>d/2+1$, and $X_S(0)=0$, $\dot{X}_S(0)=v_{S_0}\in\mathbb{R}^d$ is an arbitrary initial condition for the solid motion, then there exists a solution 
\begin{equation*}
\begin{aligned}
&(\zeta,\overline{V})\in \mathcal{C}([0,T];H^s(\mathbb{R}^d))\cap \mathcal{C}^1([0,T];H^{s-1}(\mathbb{R}^d)),\\ &X_S\in\mathcal{C}^2([0,T]),
\end{aligned}
\end{equation*}
to (\ref{eq_systemorder1}) for a sufficiently small time $T>0$ independent of $\mu$.
\end{theo}

\prove The demonstration is based on the fixed point theorem applied to an iterative scheme presented in the following subsections. The brief outline of our proof is as follows:
\begin{enumerate}
\item Reformulation of the system,
\item Construction of the iterative scheme,
\item Existence and a priori estimates for the iterative scheme,
\item Convergence of the iterative scheme solutions.
\end{enumerate}

\subsubsection{Reformulation of the coupled fluid-solid system}\label{sec_saintvenantrewrite}

Let us remark the following: the nonlinear Saint-Venant equations (\ref{eq_nonlinearSV}) admit a quasilinear hyperbolic structure. More precisely, we have the following classical reformulation using the new variable $\mathcal{U} = (\zeta,\, \overline{V})^{\top}\in\mathbb{R}^{d+1}$:
\begin{equation}\label{eq_generalhyperbolic}
\partial_t\mathcal{U} + \sum_{j=1}^{d}A_j(\mathcal{U},X_S)\partial_j\mathcal{U}+B(\mathcal{U},X_S)=0.
\end{equation}
Let us take the following real valued $(d+1)\times(d+1)$ matrices
\begin{equation}\label{matrixofnonlinearity}
A_j(\mathcal{U},X_S) = 
\left(\begin{array}{c|ccc}
\overline{V}_j & & hI_j & \\
\hline
 & & & \\
I_j^{\top} & & \overline{V}_j\operatorname{Id}_{d\times d} & \\
 & & &
\end{array}\right) \textnormal{ for } 1\leqslant j\leqslant d,
\end{equation}
where for every $1\leqslant j\leqslant d$ we have $I_j=e_j\in\mathbb{R}^d$ the $j$\textsuperscript{th} coordinate vector with respect to the standard Euclidean basis of $\mathbb{R}^d$.

We recall that $h=1+\zeta-b$ thus implying that the matrices $A_j(\mathcal{U},X_S)$ indeed depend on $X_S$, however only through the bottom variation (\ref{form_ndbottomvari}).

Following the notation in (\ref{eq_generalhyperbolic}), the additional term $B(\mathcal{U},X_S)$ is the vector 
\begin{equation*}
B(\mathcal{U},X_S)=\left(-\overline{V}\cdot\nabla_x\mathfrak{b}(x-X_S)+\nabla_x\mathfrak{b}(x-X_S)\cdot\dot{X}_S,0,\ldots,0\right)^{\top}.
\end{equation*}

From here on, we shall also use the following uniform notation for the coordinate functions of $\mathcal{U}$:
\begin{equation*}
\mathcal{U}_0=\zeta,\quad \mathcal{U}_j=\overline{V}_j\textnormal{ for }1\leqslant j\leqslant d.
\end{equation*}

As for the initial values, we have $\mathcal{U}(0,\cdot)=\mathcal{U}_{in}=(\zeta_{in},\,\overline{V}_{in})$ and $X_S(0)=0$, $\dot{X}_S(0)=v_{S_0}$. There is no restriction necessary on the initial values concerning the solid motion.

There exists a symmetrizer matrix $S(\mathcal{U},X_S)$ defined by
\begin{equation}\label{remark_reform1}
S(\mathcal{U}, X_S) = 
\left(\begin{array}{c|ccc}
1 & & 0 & \\
\hline
 & & & \\
0 & & h\operatorname{Id}_{d\times d} & \\
 & & &
\end{array}\right),
\end{equation}
such that the matrices $S(\mathcal{U},X_S)A_j(\mathcal{U},X_S)$ are symmetric. Moreover, based on our imposed lower boundary condition on $h_{in}$, one can establish that $S(\mathcal{U}_{in},0)\geqslant \min(1,h_{min})\operatorname{Id}_{(d+1)\times(d+1)}$, which guarantees that the matrix is positive definite.

Owing to the existence of such a symmetrizer $S$, the local well-posedness for a bottom with a prescribed motion follows from classical results \cite{taylorpde}. In our case and additional step is needed due to the presence of the coupling with the equation describing the solid motion.

Let us make one further remark, concerning the second order (nonlinear) ordinary differential equation characterizing the displacement of the solid $X_S$ in (\ref{eq_systemorder1o2}). Let us define the functional $\mathcal{F}[\mathcal{U}](t,Y,Z)$ as
\begin{equation*}
\mathcal{F}[\mathcal{U}](t,Y,Z)=-\frac{c_{fric}}{\sqrt{\mu}}\left(c_{solid}+\frac{1}{\tilde{M}}\int\limits_{\operatorname{supp}(\mathfrak{b})+ Y}\mathcal{U}_0\,dx\right)\frac{Z}{\left|Z\right|+\overline{\delta}}+\frac{1}{\tilde{M}}\int_{\mathbb{R}^d}\mathcal{U}_0\nabla_x\mathfrak{b}(x-Y)\,dx.
\end{equation*}

\textit{The coupled system (\ref{eq_systemorder1}) has the following equivalent form}
\begin{subequations}\label{eq_systemuniform1}
\begin{align}[left = {\,\empheqlbrace}]
&\partial_t\mathcal{U} + \sum_{j=1}^{d}A_j(\mathcal{U},X_S)\partial_j\mathcal{U}+B(\mathcal{U},X_S)=0, \label{eq_systemuniform1o1} \\
&\ddot{X}_S = \mathcal{F}[\mathcal{U}]\left(t,X_S,\dot{X}_S\right).\label{eq_systemuniform1o2}
\end{align}
\end{subequations}

\subsubsection{The iterative scheme}

To solve the coupled system (\ref{eq_systemuniform1}) we construct a sequence $\left(\{\mathcal{U}^k(t,x)\},\{X^k(t)\}\right)_{k\in\mathbb{N}}$ of approximate solutions through the scheme
\begin{subequations}\label{eq_iterativesystem1}
\begin{align}[left = {\, \empheqlbrace}]
&\;S(\mathcal{U}^k,X^k)\partial_t\mathcal{U}^{k+1} + \sum_{j=1}^{d}S(\mathcal{U}^k,X^k)A_j(\mathcal{U}^{k},X^k)\partial_j\mathcal{U}^{k+1} =-S(\mathcal{U}^k,X^k)B(\mathcal{U}^k,X^k), \label{eq_iterativesystem1o1} \\
&\;\ddot{X}^{k+1} = \mathcal{F}[\mathcal{U}^{k+1}]\left(t,X^{k+1},\dot{X}^{k+1}\right); \label{eq_iterativesystem1o2} \\
&\;\mathcal{U}^{k+1}(0,\cdot)=\mathcal{U}_{in},\quad X^{k+1}(0)=0,\;\dot{X}^{k+1}(0)=v_{S_0}. \nonumber
\end{align}
\end{subequations}
Here the matrices $A_j$ and $S$ are the matrices defined in (\ref{matrixofnonlinearity}) and (\ref{remark_reform1}). In what follows we will make use of the following abbreviations
\begin{equation*}
S^k=S(\mathcal{U}^k,X^k),\quad A_j^k=A_j(\mathcal{U}^{k},X^k),\;\textnormal{ and }\;B^k=B(\mathcal{U}^k,X^k).
\end{equation*}

The main goal is to prove the existence and convergence of this sequence. We will follow the footsteps of a classical method, presented by Alinhac and Gérard in \cite{alinhacgerard} for instance, detailing only the parts where additional estimates are necessary due to the coupling terms.

The iterative scheme works as follows: we choose the initial $k=0$ elements to be $(\mathcal{U}^0,X^0)=(\mathcal{U}_{in},0)$. From then on, at each step $k$ ($k\in\mathbb{N}$) we have to solve a linear symmetric hyperbolic PDE system (\ref{eq_iterativesystem1o1}) to recover $\mathcal{U}^{k+1}$, and then a second order nonlinear ODE (\ref{eq_iterativesystem1o2}) to obtain $X^{k+1}$.

\subsubsection{Existence and a priori estimates}

Now, the aim is to establish the existence of solutions $(\mathcal{U}^{k+1},X^{k+1})$ ($k\geqslant 0$) for the iterative scheme to justify their definition in (\ref{eq_iterativesystem1}). Furthermore we shall also obtain a control of the velocity fields for our coupled system. In particular an upper bound on $\mathcal{U}^{k+1}$ in a ``large norm'', partially in order to guarantee the boundedness conditions required for the existence result presented, as well as to introduce certain inequalities which will be useful for the convergence of the series.

In what follows, we will make use of the following notation
\begin{defi}\label{def_bracket}
For an $f(t,x)\in L^\infty([0,T];\mathcal{X}(\mathbb{R}^d))$ function let us define
\begin{equation}
\brck{f}_{T,H^s}=\sup_{t\in[0,T]}\|f(t,\cdot)\|_{\mathcal{X}}.
\end{equation}
\end{defi}

With this definition at our disposal, we can state the induction hypothesis (\ref{recurrenthyp}) for the boundedness of solutions $(\mathcal{U}^l,X^l)_{l\leqslant k}$ of (\ref{eq_iterativesystem1}):
\begin{equation}\label{recurrenthyp}\tag{$H_k$}
\textnormal{for } 0\leqslant l\leqslant k,\quad \brck{\mathcal{U}^l}_{T,H^s}\leqslant C_f,\;\brck{\mathcal{U}^l-\mathcal{U}_{in}}_{T,L^\infty}\leqslant\delta_0,\; \sup_{t\in[0,T]}\left|\dot{X}^l-v_{S_0}\right|\leqslant C_vT,
\end{equation}
for a sufficiently large constant $C_f=C(\tilde{M}^{-1},S_0^{-1},\|\mathcal{U}_{in}\|_{H^s},\|\mathfrak{b}\|_{H^s})$, with $\delta_0>0$ a small constant to be defined, independent of $k$, and $C_v=C(\tilde{M}^{-1},\|\mathfrak{b}\|_{H^s};C_f)$,

\begin{prop}\label{stat_recurrenceSV}
For $k\geq 0$, assuming (\ref{recurrenthyp}), there exists a solution $\mathcal{U}^{k+1}\in\mathcal{C}([0,T];H^s(\mathbb{R}^d))$, $X^{k+1}\in\mathcal{C}^2([0,T])$ of (\ref{eq_iterativesystem1}), moreover, by an adequate choice of $C_f$, $\delta_0$, and $T$ (independent of $k$)
\begin{equation*}
(H_k)\Rightarrow(H_{k+1}).
\end{equation*}
\end{prop}

\prove The proof goes by induction. For $k=0$, $(H_0)$ is clearly verified. For the induction step, we shall treat separately the case of the PDE (part A) and the case of the ODE (part B), for the sake of clarity.

\noindent\textbf{Part A: existence and energy estimate for $\mathcal{U}^{k+1}$:} The initial values $\mathcal{U}^{k+1}(0,\cdot)$ are bounded since they are equal to the original initial values $\mathcal{U}_{in}$. Since we are operating by induction with respect to $k$, for the respective $\mathcal{U}^k$ term we already have existence, moreover we also have the large norm estimates (\ref{recurrenthyp}) at hand, which in particular guarantees the uniform bounds for $\mathcal{U}^k$ (independently of the index $k$) for small time $T$ and $\delta_0$. Also, given the simple structure of $S^k$ and $S^kA_j^k$, they are bounded as well in Lipschitz norm.

\begin{lem}
For $k\geq 0$, with the initial condition $\mathcal{U}^{k+1}(0,\cdot)=\mathcal{U}_{in}$ and the hypothesis (\ref{recurrenthyp}) there exists a $\mathcal{C}([0,T];H^s(\mathbb{R}^d))$ solution $\mathcal{U}^{k+1}$ for the linear symmetric hyperbolic PDE system defined in (\ref{eq_iterativesystem1o1}).
\end{lem}
\prove Notice that (\ref{eq_iterativesystem1o1}) has a particular symmetric structure which may be exploited based on the following proposition:
\begin{prop}\label{prop_structurelinear}
Let us consider the symmetric hyperbolic differential operator 
\begin{equation*}
L=S\pt+\sum_{j=1}^dSA_j\partial_j
\end{equation*}
with $S$ and $SA_j$ symmetric real valued and bounded in Lipschitz norm, with $S\geqslant S_0\operatorname{Id}$, where $S_0>0$ over $[0,T]$. Furthermore let us consider $s\in\mathbb{R}$, $s>d/2+1$ and let us take
\begin{equation*}
\lambda_s=C\left(\brck{S}_{T,H^s},\brck{\sum_{j=1}^dSA_j}_{T,H^s},\brck{\pt S}_{T,L^\infty}\right).
\end{equation*}
Then, for any $f\in L^1([0,T];H^s(\mathbb{R}^d))$ and $\varphi\in H^s(\mathbb{R}^d)$ the Cauchy problem
\begin{equation}\label{eq_linhypsym}
\begin{cases}
Lu=f,\;0<t<T\\
u(0,\cdot)=\varphi,
\end{cases}
\end{equation}
admits a unique solution $u\in\mathcal{C}([0,T];H^s(\mathbb{R}^d))$ that verifies the energy estimate
\begin{equation}\label{eq_enerestimlin}
S_0\sup_{t\in[0,T]}\left\{e^{-\lambda_s t}\|u(t,\cdot)\|_{H^s}\right\}\leqslant\|\varphi\|_{H^s}+2\int_0^Te^{-\lambda_s t'}\|f(t',\cdot)\|_{H^s}dt'.
\end{equation}
\end{prop}
For more details as well as a complete proof, we refer to \cite{alinhacgerard}.
\begin{rem}
Under the same regularity assumptions, we can also infer that
\begin{equation}\label{eq_enerestimlin0}
S_0\sup_{t\in[0,T]}\left\{e^{-\lambda_0 t}\|u(t,\cdot)\|_{L^2}\right\}\leqslant\|\varphi\|_{L^2}+2\int_0^Te^{-\lambda_0 t'}\|f(t',\cdot)\|_{L^2}dt',
\end{equation}
where
\begin{equation*}
\lambda_0=\frac{1}{2S_0}\brck{\frac{1}{2}\sum_{j=1}^d\partial_j(SA_j)-\partial_t S}_{T,L^\infty}.
\end{equation*}
\end{rem}

We want to apply Proposition \ref{prop_structurelinear}. to solve the linear PDE (\ref{eq_iterativesystem1o1}) for $\mathcal{U}^{k+1}$ in $\mathcal{C}([0,T];H^s(\mathbb{R}^d))$. First of all, we have that $\mathcal{U}^k$ and $X^k$ are continuous in time. So, since $h_{in}\geqslant h_{min}$, by using (\ref{recurrenthyp}) for a $\delta_0$ sufficiently small, we obtain $S^k\geqslant S_0\operatorname{Id}$ in $[0,T]$ for $S_0=(1/2)\min(1,h_{min})$.

By the regularity of $\mathcal{U}^k$ and $\mathfrak{b}$, the source term $S^kB^k$ is also in $H^s$. More exactly, we have the following
\begin{lem}\label{lemma_rhsestimate}
The source term of (\ref{eq_iterativesystem1o1}) satisfies the following linear-in-time estimate
\begin{equation}
\|S^kB^k\|_{H^s}\leqslant C_F(1+T),
\end{equation}
with the constant $C_F=C(\|\mathfrak{b}\|_{H^s};C_f,C_v)$, independent of $k$ and of $T$.
\end{lem}
\prove We have
\begin{equation*}
\begin{aligned}
\|S^kB^k\|_{H^s}&\lesssim\left(1+\|\mathcal{U}^k_0\|_{L^{\infty}}+\|\mathfrak{b}\|_{L^\infty}\right)\cdot\|B^k\|_{H^s}+\left(1+\|\mathcal{U}^k_0\|_{H^s}+\|\mathfrak{b}\|_{H^s}\right) \cdot\|B^k\|_{L^\infty}\\
&\lesssim \left(1+\|\mathcal{U}^k_0\|_{H^s}+\|\mathfrak{b}\|_{H^s}\right)\cdot\|\mathfrak{b}\|_{H^s}\left(\|\overline{V}^k\|_{H^s}+\left|\dot{X}^k\right|\right)\\
\end{aligned}
\end{equation*}
using the special structure of the matrix $S^k$, the Sobolev embedding $H^s(\mathbb{R}^d)\hookrightarrow L^\infty(\mathbb{R}^d)$ (which is valid since $s>d/2$) as well as the fact that the Sobolev norm is translation invariant. Then, the induction hypothesis (\ref{recurrenthyp}) provides a uniform bound $C_f$ for $\|\mathcal{U}^k\|_{H^s}$, as well as a linear-in-time estimate for $|\dot{X}^k|$, so, since $\mathfrak{b}$ is still regular,
\begin{equation*}
\|S^kB^k\|_{H^s}\leqslant C(\|\mathfrak{b}\|_{H^s};C_f, C_v)(1+T).
\end{equation*} \quod

Now we only need to verify that $\lambda_s$ is bounded. For this, we have that
\begin{lem}
Assuming that (\ref{recurrenthyp}) holds,
\begin{equation*}
\lambda_s\leqslant c(\|\mathfrak{b}\|_{H^s};C_f, C_v)(1+T),
\end{equation*}
where $c$ is a continuous nondecreasing function of its arguments.
\end{lem}
\prove Making use of the fact that $S^k$ and $A_j^k$ depend on $X^k$ in a very simple way, throughout the function $h^k$, thus it is present as a translation for the function $\mathfrak{b}$, which obviously does not affect the $L^{\infty}$ or $H^s$ norms, we have that
\begin{equation*}
\brck{S}_{T,H^s}\leqslant 1+\brck{\mathcal{U}}_{T,H^s}+\|\mathfrak{b}\|_{H^s},
\end{equation*}
and that
\begin{equation*}
\|SA_j\|_{H^s}\leqslant \|S\|_{L^\infty}\|A_j\|_{H^s}+\|S\|_{H^s}\|A_j\|_{L^\infty}\leqslant c(\|\mathfrak{b}\|_{H^s}, \|\mathcal{U}\|_{H^s}),
\end{equation*}
by the Sobolev embedding $H^s(\mathbb{R}^d)\hookrightarrow L^{\infty}(\mathbb{R}^d)$.

As for the estimate on $\|\pt S^k\|_{L^\infty}$, we estimate the $L^\infty$ norm of $\pt h^k$, which is
\begin{equation*}
\pt h^k = \pt\mathcal{U}_0^k+\nabla_x\mathfrak{b}(x-X^k)\cdot\dot{X}^k.
\end{equation*}

The second term is already controlled by (\ref{recurrenthyp}). Based on the corresponding equation for $\mathcal{U}^k$ (from equation (\ref{eq_iterativesystem1o1})), we have that 
\begin{equation*}
\pt\mathcal{U}^k+\sum_{j=1}^{d}A_j^{k-1}\partial_j\mathcal{U}^k=-B^{k-1},
\end{equation*}
which implies that
\begin{equation*}
\begin{aligned}
\|\pt\mathcal{U}^k\|_{L^\infty}&\leqslant \sum_{j=1}^d\|A_j^{k-1}\partial_j\mathcal{U}^k\|_{L^\infty}+\left\|\overline{V}^{k-1}\cdot\nabla_x\mathfrak{b}(.-X^{k-1})+ \nabla_x\mathfrak{b}(.-X^{k-1})\cdot\dot{X}^{k-1}\right\|_{L^\infty}\\
&\lesssim \sum_{j=1}^d\left(1+\|\mathcal{U}^{k-1}\|_{L^{\infty}}+\|\mathfrak{b}\|_{L^{\infty}}\right)\cdot\|\mathcal{U}^k\|_{H^s}+\|\mathfrak{b}\|_{H^s}\left(\|\mathcal{U}^{k-1}\|_{L^\infty}+\left|\dot{X}^{k-1}\right|\right)\\
&\leqslant c(\|\mathfrak{b}\|_{H^s};C_f, C_v)(1+T),
\end{aligned}
\end{equation*}
by (\ref{recurrenthyp}) and the regularity of $\mathfrak{b}$, as well as the Sobolev embedding $H^s(\mathbb{R}^d)\hookrightarrow W^{1,\infty}(\mathbb{R}^d)$ ($s>d/2+1$). Therefore, $\lambda_s$ is indeed a constant independent of $k$, and linear in $T$. \quod

Now, we turn our attention towards the first two estimates in $(H_{k+1})$. For the large norm estimate, the $H^s$ energy estimate (\ref{eq_enerestimlin}) from Proposition \ref{prop_structurelinear}. for equation (\ref{eq_iterativesystem1o1}) of $\mathcal{U}^{k+1}$ can be stated to obtain
\begin{equation}\label{eq_estimationresidual}
S_0\sup_{t\in[0,T]}\left\{e^{-\lambda_s t}\|\mathcal{U}^{k+1}(t,\cdot)\|_{H^s}\right\}\leqslant \|\mathcal{U}_{in}\|_{H^s} + 2\int_0^Te^{-\lambda_s t}\|S^kB^k\|_{H^s}dt.
\end{equation}

The right hand side of (\ref{eq_iterativesystem1o1}) can be estimated by Lemma \ref{lemma_rhsestimate}, so we obtain
\begin{equation*}
\brck{\mathcal{U}^{k+1}}_{T,H^s}\leqslant \frac{1}{S_0}e^{\lambda_s T}\|\mathcal{U}_{in}\|_{H^s}+T(1+T) e^{\lambda_s T} c(\|\mathfrak{b}\|_{H^s};C_f, C_v).
\end{equation*}
For $C_f$ sufficiently large the first term in the right hand side is less than $C_f/2$. Therefore, for $T$ small enough the second term will be less than $C_f/2$ too. This proves the first estimate of $(H_{k+1})$.

In order to obtain a uniform $L^\infty$ estimate for $\mathcal{U}^{k+1}-\mathcal{U}_{in}$, we shall first of all control $\pt \mathcal{U}^{k+1}$ in $L^\infty$. Just as before for $\|\pt\mathcal{U}^k\|_{L^\infty}$, by the large norm estimate for $\mathcal{U}^{k+1}$ we have that
\begin{equation*}
\|\pt\mathcal{U}^{k+1}\|_{L^\infty}\leqslant c(\|\mathfrak{b}\|_{H^s};C_f, C_v)(1+T).
\end{equation*}
Therefore we obtain
\begin{equation*}
\brck{\mathcal{U}^{k+1}-\mathcal{U}_{in}}_{T,L^\infty}\leqslant T\|\pt\mathcal{U}^{k+1}\|_{L^\infty}\leqslant c(\|\mathfrak{b}\|_{H^s};C_f, C_v)T(1+T).
\end{equation*}
Hence, for a sufficiently small time $T$ (independently of $k$) we get that the right hand side is less than $\delta_0$.

\noindent\textbf{Part B: existence and velocity estimate for $X^{k+1}$:} For the existence of the solution $X^{k+1}$ of the ODE (\ref{eq_iterativesystem1o2}), we shall apply the Picard--Lindel\"{o}f theorem.
\begin{lem}
For $k\geq 0$, with the initial conditions $X^{k+1}(0)=0$, and $\dot{X}^{k+1}(0)=v_{S_0}$ and the hypothesis (\ref{recurrenthyp}), there exists a continuously differentiable solution $X^{k+1}$ for the nonlinear second order non-homogeneous ODE defined in (\ref{eq_iterativesystem1o2}) for $t\in[0,T_S]$, where $T_S=C(\|\mathfrak{b}\|_{H^s};C_f)$.
\end{lem}
\prove For the Picard-Lindel\"{o}f theorem, we have to show that the nonlinear functional on the right hand side of (\ref{eq_iterativesystem1o2}) is continuous in time and uniformly Lipschitz in the spatial variable.

We recall that the functional $\mathcal{F}[\mathcal{U}^{k+1}](t,Y,Z)$ has the form of
\begin{equation*}
\mathcal{F}[\mathcal{U}^{k+1}](t,Y,Z)=-\frac{c_{fric}}{\sqrt{\mu}}\left(c_{solid}+\frac{1}{\tilde{M}}\int\limits_{\operatorname{supp}(\mathfrak{b})+ Y}\mathcal{U}_0^{k+1}\,dx\right)\frac{Z}{\left|Z\right|+\overline{\delta}}+\frac{1}{\tilde{M}}\int_{\mathbb{R}^d}\mathcal{U}_0^{k+1}\nabla_x\mathfrak{b}(x-Y)\,dx.
\end{equation*}

We already know that $\mathcal{U}^{k+1}_0$ is of class $\mathcal{C}([0,T];H^s(\mathbb{R}^d)$, so it is continuous in the time variable, regular in the spatial variable, moreover the function $\mathfrak{b}$ is regular and with a compact support, thus the integrals indeed exist and are bounded, furthermore based on the well known theorem concerning the continuity of a parametric integral, it will be continuous with respect to $t$.

All we need to show is that it is (locally) uniformly Lipschitz with respect to its second variable $(Y,Z)$. Examining $\mathcal{F}[\mathcal{U}^{k+1}](t,Y,Z)$, it is clear that the second term is Lipschitz continuous due to the fact that $\mathfrak{b}$ is regular. As for the first term, since it contains a product of multiple terms with the variables $Y$ and $Z$, by adding and subtracting intermediate terms they can be separated.

Let us take a closer look on these two separate terms. The integral term can be estimated due to
\begin{equation*}
\int\limits_{\operatorname{supp}(\mathfrak{b})+ Y}\mathcal{U}_0^{k+1}(t,x)\,dx=\int\limits_{\operatorname{supp}(\mathfrak{b})}\mathcal{U}_0^{k+1}(t,x- Y)\,dx,
\end{equation*}
and the regularity of $\mathcal{U}^{k+1}_0$. Since we chose $\overline{\delta}>0$, the function 
\begin{equation*}
Z\mapsto\frac{Z}{|Z|+\overline{\delta}}
\end{equation*}
is Lipschitz continuous. So, putting all the estimates together, we obtain that
\begin{equation*}
\begin{aligned}
\bigg|\mathcal{F}[\mathcal{U}^{k+1}](t,Y_1,Z_1)-\mathcal{F}[\mathcal{U}^{k+1}](t,Y_2,Z_2)\bigg|&\leqslant \frac{3c_{fric}}{\overline{\delta}\sqrt{\mu}}\left(c_{solid}+\frac{C_f|\operatorname{supp}(\mathfrak{b})|}{\tilde{M}}\right)\left|Z_1-Z_2\right|\\
&\quad\,+2C_fL_{\nabla_x\mathfrak{b}}\frac{|\operatorname{supp}(\mathfrak{b})|}{\tilde{M}}\cdot|Y_1-Y_2|;
\end{aligned}
\end{equation*}
where $L$ denotes the Lipschitz constant of the corresponding function in the subscript.\quod

One of the most important parts of the proof is the control on the solid velocity, since the solid equation contains an order $\mu^{-1/2}$ term which could potentially become huge, making the system blow up. However,
\begin{lem}\label{lemma_velocityestimate}
A control on the solid velocity is ensured by
\begin{equation}
\left|\dot{X}^{k+1}(t)\right|\leqslant \left|v_{S_0}\right|+C_vt,
\end{equation}
with a constant $C_v=C(\tilde{M}^{-1},\|\mathfrak{b}\|_{H^s};C_f)$ independent of $k$, $\mu$ and $t$.
\end{lem}
\prove By definition $X^{k+1}$ satisfies the corresponding second order nonlinear non-homogeneous equation in (\ref{eq_iterativesystem1o2}) so we have that
\begin{equation*}
\ddot{X}^{k+1}=\mathcal{F}[\mathcal{U}^{k+1}]\left(t,X^{k+1},\dot{X}^{k+1}\right),
\end{equation*}
thus, multiplying by $\dot{X}^{k+1}(t)$, we get that
\begin{equation*}
\begin{aligned}
\ddot{X}^{k+1}\cdot\dot{X}^{k+1} &= -\frac{c_{fric}}{\sqrt{\mu}}\left(c_{solid}+\frac{1}{\tilde{M}}\int\limits_{\operatorname{supp}(\mathfrak{b})+X^{k+1}}\mathcal{U}_0^{k+1}\,dx\right)\frac{\left|\dot{X}^{k+1}\right|^2}{\left|\dot{X}^{k+1}\right|+\overline{\delta}}+\frac{1}{\tilde{M}}\int_{\mathbb{R}^d}\mathcal{U}^{k+1}_0\nabla_x\mathfrak{b}(x-X^{k+1})\,dx\cdot\dot{X}^{k+1}\\
&\leqslant 0+\frac{1}{\tilde{M}}\int_{\mathbb{R}^d}\mathcal{U}^{k+1}_0\nabla_x\mathfrak{b}(x-X^{k+1})\cdot\dot{X}^{k+1}\,dx,
\end{aligned}
\end{equation*}
here the key remark is that the first, negative term disappeared from the equations. Thus, we are left with
\begin{equation*}
\frac{1}{2}\frac{d}{dt}\left[\left|\dot{X}^{k+1}\right|^2\right]\leqslant\frac{1}{\tilde{M}}\|\mathcal{U}^{k+1}_0\|_{L^2}\|\nabla_x\mathfrak{b}\|_{L^2}\cdot|\dot{X}^{k+1}|\leqslant\frac{C_f}{\tilde{M}}\|\mathfrak{b}\|_{H^s}\left|\dot{X}^{k+1}\right|,
\end{equation*}
so by a Gr\"onwall type lemma for $\dot{X}^{k}$, we may conclude that
\begin{equation*}
\left|\dot{X}^{k}(t)\right|\leqslant |v_{S_0}| + \frac{C_f}{\tilde{M}}\|\mathfrak{b}\|_{H^s}t.
\end{equation*}

\quod

This concludes the velocity estimate for the object. So we proved the existence of solutions $(\mathcal{U}^{k+1},X^{k+1})$ for the system (\ref{eq_iterativesystem1}), moreover we established the necessary elements for the upper bounds concerning the velocities in $(H_{k+1})$. \quod

\subsubsection{Convergence}

We want to establish the convergence of the series from (\ref{eq_iterativesystem1}), for that we need the $L^2$-norm estimates for the difference between two subsequent elements for $\mathcal{U}^k$, for $X^k$ we shall simply estimate in $\mathbb{R}^d$-norm.

We start by subtracting the equations corresponding to the $k$\textsuperscript{th} element from the equations corresponding to the $(k+1)$\textsuperscript{th} element. After the subtraction we have
\begin{subequations}\label{eq_convergence1}
\begin{align}[left = {\, \empheqlbrace}]
\begin{split}
S^k\pt(\mathcal{U}^{k+1}&-\mathcal{U}^k) + \sum_{j=1}^dS^kA_j^k\partial_j(\mathcal{U}^{k+1}-\mathcal{U}^k) =  \\
= &-(S^k-S^{k-1})\pt\mathcal{U}^k-\sum_{j=1}^d(S^kA_j^k-S^{k-1}A_j^{k-1})\partial_j\mathcal{U}^k - (S^kB^k-S^{k-1}B^{k-1}), 
\end{split}\label{eq_convpde1}\\
(\mathcal{U}^{k+1}-\mathcal{U}^k)(0,\cdot) &=0; \nonumber \\
\frac{d^2}{dt^2}(X^{k+1}-X^k)&= \mathcal{F}[\mathcal{U}^{k+1}]\left(t,X^{k+1},\dot{X}^{k+1}\right)-\mathcal{F}[\mathcal{U}^{k}]\left(t,X^{k},\dot{X}^{k}\right) \label{eq_convode1} \\
(X^{k+1}-X^k)(0)&= 0,\;\quad\frac{d}{dt}(X^{k+1}-X^k)(0)=0. \nonumber
\end{align}
\end{subequations}

We provide separately an appropriate estimate in this small norm for the solutions $(\mathcal{U}^{k+1}-\mathcal{U}^k)$ and $(X^{k+1}-X^k)$ of the system. The estimate for the ODE part (\ref{eq_convode1}) is given by the following lemma.
\begin{lem}
For a solution $X^{k+1}-X^k$ of (\ref{eq_convode1}), we have that
\begin{equation}
\sup_{t\in[0,T]}|X^{k+1}(t)- X^k(t)|\lesssim C_1T\sup_{t'\in[0,T]}\left(\left\|\mathcal{U}^{k+1}(t',\cdot)-\mathcal{U}^{k}(t',\cdot)\right\|_{L^2}+\left|X^{k+1}(t')-X^{k}(t')\right|\right),
\end{equation}
with a constant $C_1=C(\tilde{M}^{-1},\|\mathfrak{b}\|_{H^s};C_f)$.
\end{lem}

Notice that the right hand side of the estimate contains exactly the same differences as the ones we would like to establish an upper bound for, but due to the presence of the factor $T$, with $T$ sufficiently small, it will be completely absorbed by the left hand side.

\prove To treat the difference of products that arise multiple times, we introduce intermediary terms, just as we did for the verification of the Lipschitz-property. So following standard computations, we get that
\begin{equation*}
\begin{aligned}
|X^{k+1}(t)&- X^k(t)|=\left|\int_0^t\frac{d}{dt}(X^{k+1}-X^k)(s)\,ds\right|=\left|\int_0^t\int_0^s\frac{d^2}{dt^2}(X^{k+1}-X^k)(\tau)\,d\tau ds\right|\\
&\lesssim \int_0^t\int_0^s(1+T)\|\mathcal{U}^{k+1}\|_{L^2}|X^{k+1}-X^{k}|(\tau)\,d\tau ds+\int_0^t\left|\int_0^s\frac{d}{d\tau}(X^{k+1}-X^{k})(\tau)\,d\tau\right| ds\\
&\,\quad+ \int_0^t\int_0^s\sup_{t'\in[0,T]}\left\|\mathcal{U}^{k+1}_0(t',\cdot)-\mathcal{U}^{k}_0(t',\cdot)\right\|_{L^2}\|\nabla_x\mathfrak{b}\|_{L^2}\,d\tau ds+\int_0^t\int_0^s |X^{k+1}-X^{k}|(\tau)\,d\tau ds \\
\end{aligned}
\end{equation*}
Here we estimate each term by the controls of (\ref{recurrenthyp}) and the Lipschitz properties, in order to obtain
\begin{equation*}
|X^{k+1}(t)- X^k(t)|\lesssim T\left(\sup_{t'\in[0,T]}\left\|\mathcal{U}^{k+1}(t',\cdot)-\mathcal{U}^{k}(t',\cdot)\right\|_{L^2}+\sup_{t'\in[0,T]}|X^{k+1}(t')-X^{k}(t')|\right).
\end{equation*}
\quod

Then, for the estimate for the PDE part we have the following.
\begin{lem}
For a solution $ $ of (\ref{eq_convpde1}), we obtain
\begin{equation}
\sup_{t\in[0,T]}\big\|\mathcal{U}^{k+1}(t,\cdot) - \mathcal{U}^k(t,\cdot)\big\|_{L^2}\lesssim C_2T\sup_{t'\in[0,T]}\left(\left\|\mathcal{U}^k(t',\cdot)-\mathcal{U}^{k-1}(t',\cdot)\right\|_{L^2}+\left|X^k(t')-X^{k-1}(t')\right|\right),
\end{equation}
with a constant $C_2=C(\tilde{M}^{-1},\|\mathfrak{b}\|_{H^s};C_f, C_v)$.
\end{lem}
\prove Once again we aim to use the energy estimate, since we have the linear system (\ref{eq_convpde1}) of the type of Proposition \ref{prop_structurelinear}. for the variable $\mathcal{U}^{k+1}-\mathcal{U}^k$. 

The same reasoning applies here as for the energy estimates section concerning the applicability of the proposition, since $\lambda_0$ is bounded, so we only need a sufficient upper bound for the right hand side of (\ref{eq_convpde1}), denoted by $F$. We shall examine it term by term.

The first two terms are handled with standard techniques to deduce
\begin{equation*}
\|(S^k-S^{k-1})\pt\mathcal{U}^k\|_{L^2}\lesssim\|\mathcal{U}^k-\mathcal{U}^{k-1}\|_{L^2}+|X^k-X^{k-1}|,
\end{equation*}
\begin{equation*}
\left\|(S^kA_j^k-S^{k-1}A_j^{k-1})\partial_j\mathcal{U}^k\right\|_{L^2}\lesssim\left\|\mathcal{U}^k-\mathcal{U}^{k-1}\right\|_{L^2}+|X^k-X^{k-1}|.
\end{equation*}

And for the third term, we deal with the product via an intermediary term, so 
\begin{equation*}
\|S^kB^k-S^{k-1}B^{k-1}\|_{L^2}\leqslant\|S^k-S^{k-1}\|_{L^2}\|B^k\|_{L^\infty}+\|S^{k-1}\|_{L^\infty}\|B^k-B^{k-1}\|_{L^2},
\end{equation*}
then, by using the definition of the source term, we get that
\begin{equation*}
\begin{aligned}
\|S^kB^k-S^{k-1}B^{k-1}\|_{L^2}&\lesssim(\|\mathcal{U}^k-\mathcal{U}^{k-1}\|_{L^2}+|X^k-X^{k-1}|)\|\overline{V}^k\cdot\nabla_x\mathfrak{b}(.- X^k)\|_{L^\infty}\\
&\,\quad+(\|\mathcal{U}^k-\mathcal{U}^{k-1}\|_{L^2}+|X^k-X^{k-1}|)\left\|\nabla_x\mathfrak{b}(.- X^k)\cdot\dot{X}^k\right\|_{L^\infty}\\
&\,\quad+(1+\|\mathcal{U}^{k-1}\|_{L^\infty})\left\|\overline{V}^k\cdot\nabla_x\mathfrak{b}(.- X^k)-\overline{V}^{k-1}\cdot\nabla_x\mathfrak{b}(.- X^{k-1})\right\|_{L^2}\\
&\,\quad+(1+\|\mathcal{U}^{k-1}\|_{L^\infty})\left\|\nabla_x\mathfrak{b}(.- X^k)\cdot\dot{X}^k -\nabla_x\mathfrak{b}(.- X^{k-1})\cdot\dot{X}^{k-1}\right\|_{L^2}.\\
\end{aligned}
\end{equation*}

Again, with the apparition of intermediary terms for each product, by utilizing Lemma \ref{lemma_velocityestimate}, we get that
\begin{equation*}
\|S^kB^k-S^{k-1}B^{k-1}\|_{L^2}\lesssim (1+T)\sup_{t'\in[0,T]}\left(\|\mathcal{U}^k(t',\cdot)-\mathcal{U}^{k-1}(t',\cdot)\|_{L^2}+|X^k(t')-X^{k-1}(t')|\right).
\end{equation*}

Applying the energy estimate to (\ref{eq_convpde1}), we obtain that
\begin{equation*}
\begin{aligned}
S_0\sup_{t\in[0,T]}\Big\{e^{-\lambda t}\|\mathcal{U}^{k+1}(t,\cdot)&-\mathcal{U}^k(t,\cdot)\|_{L^2}\Big\}\leqslant 2\int_0^Te^{-\lambda t}\left\|F(t,\cdot)\right\|_{L^2}dt\\
&\lesssim T(1+T)\sup_{t'\in[0,T]}\left\{e^{-\lambda t'}\left(\left\|\mathcal{U}^k(t',\cdot)-\mathcal{U}^{k-1}(t',\cdot)\right\|_{L^2}+|X^k(t')-X^{k-1}(t')|\right)\right\}.
\end{aligned}
\end{equation*}
\quod

To sum up the results from the two previous lemmas, we obtained that for $T$ sufficiently small, we have that
\begin{equation}
\begin{aligned}
\sup_{t\in[0,T]}\Big\|\mathcal{U}^{k+1}(t,\cdot) &- \mathcal{U}^k(t,\cdot)\Big\|_{L^2}+\sup_{t\in[0,T]}|X^{k+1}(t)-X^k(t)|\\
&\leqslant\mathfrak{c}\left(\sup_{t\in[0,T]}\left\|\mathcal{U}^k(t,\cdot)-\mathcal{U}^{k-1}(t,\cdot)\right\|_{L^2}+\sup_{t\in[0,T]}|X^k(t)-X^{k-1}(t)|\right),
\end{aligned}
\end{equation}
with a constant $\mathfrak{c}<1$. This ensures the convergence in $L^{\infty}([0,T];L^2(\mathbb{R}^d))$. 

Since we have that the $H^s$-norm is bounded, we may extract a weakly convergent subsequence from the series, and since the limit in the sense of distributions is unique, we have convergence for the whole series in $H^s$ too. \quod

This concludes the proof of the theorem, since, to deduce the regularity implied in the statement, we only have to use the convexity of the norm, following classical regularity arguments. Thus we obtained a classical solution of the coupled system (\ref{eq_systemorder1}) for a sufficiently small time $T$.

\section{The $\mathcal{O}(\mu^2)$ asymptotic regime: The Boussinesq system}

In this section, we move on to the next order regarding the asymptotic regime, that is the approximations of order $\mu^2$. In order to simplify the computations, we consider here a weakly nonlinear regime, i.e. we assume that $\varepsilon=\mathcal{O}(\mu)$. The fluid is then governed by a Boussinesq system. Thus, the asymptotic regime writes as follows
\begin{equation}\tag{BOUS}\label{hyp_boussinesqregime}
0\leqslant\mu\leqslant\mu_{max}\ll 1,\quad\varepsilon=\mathcal{O}(\mu).
\end{equation}

At second order, the asymptotic expansion of $\overline{V}$ in terms of $\zeta$ and $\overline{V}$ (based on Proposition 3.37. of \cite{lannesbible}) is given by
\begin{equation*}
\overline{V}= \nabla_x\psi+\frac{\mu}{3h}\nabla_x(h^3\nabla_x\cdot\nabla_x\psi)-\frac{\mu}{2}h\nabla_x\pt b+\mathcal{O}(\mu^2),
\end{equation*} 
so by making use of the definition of $h$, once again taking the gradient of the second equation in (\ref{eq_ndzakharovmovingbottom}), and neglecting terms of order $\mathcal{O}(\mu^2)$, equations (\ref{eq_ndzakharovmovingbottom}) under the Boussinesq regime (\ref{hyp_boussinesqregime}) take the form of
\begin{equation}\label{eq_boussinesq}
\begin{cases}
\partial_t\zeta+\nabla_x\cdot(h\overline{V}) = \pt b,\\
\left(1- \displaystyle\frac{\mu}{3}\Delta_x\right)\partial_t\overline{V}+\nabla_x\zeta+\varepsilon(\overline{V}\cdot\nabla_x)\overline{V} = -\displaystyle\frac{\mu}{2}\nabla_x\partial_t^2b.
\end{cases}
\end{equation}
For the well-posedness of this Boussinesq system, see for instance \cite{lannesbible}.

\begin{rem}
Without the smallness assumption on $\varepsilon=\mathcal{O}(\mu)$, it is still possible to perform an asymptotic expansion at $\mathcal{O}(\mu^2)$. The resulting system is more general than the Boussinesq system (\ref{eq_boussinesq}) but also more complicated, it is known as the Serre--Green--Naghdi equations. For the justification of this general system in the fixed bottom case, please refer to \cite{alvarezlanneslong}, or to \cite{iguchiGNtsunami} for a moving bottom under a forced motion.
\end{rem}

It is well-known for the fixed bottom case that the good timescale of Boussinesq-type systems is of order $\varepsilon^{-1}$ in order to be able to properly observe the nonlinear and dispersive effects of equations (\ref{eq_boussinesq}) (see for instance \cite{davidboussinesq,sautxulong,cosminboussinesq}). However, for a time dependent bottom (as it is in our case), one can only infer an existence time of $\mathcal{O}(1)$, due to the source term $\partial_t b$ on the right hand side of the first equation in (\ref{eq_boussinesq}). Throughout this section, we show that, with the presence of the solid in the system as well as with better estimates, a time of existence in $\varepsilon^{-1/2}$ is achievable.

\subsection{Formal derivation of the corresponding solid motion equation}

As we may observe from the Boussinesq system (\ref{eq_boussinesq}) the bottom related source terms are respectively of order $\mathcal{O}(1)$ and $\mathcal{O}(\mu)$ for the first and second equations. To ensure at least a reasonable level of consistency on the whole coupled system, we have to impose (at least) the same precision in deriving formally the equation dealing with the solid; the surface integral present in (\ref{eq_ndnewton}) will therefore be approximated at order $\mathcal{O}(\mu^2)$.

Our strategy is exactly the same as for the first order approximation case in the previous section, but it is carried out to the next order of approximation. However, it turns out that due to the additional hypothesis on $\varepsilon$, the pressure formula (\ref{form_pressureassymptotic}) derived in Section 2 still holds in this regime, namely

\begin{lem}\label{lemma_pressureasymptotic2} Under the Boussinesq hypotheses (\ref{hyp_boussinesqregime}), the pressure takes the following form
\begin{equation}
P=\frac{P_{atm}}{\varrho gH_0}+(\varepsilon\zeta-z)+\mathcal{O}(\mu^2).
\end{equation}
\end{lem}

\prove The residual in (\ref{form_pressureassymptotic}) is of size $\mathcal{O}(\mu)$ in the \ref{hyp_saintvenantregime} Saint-Venant regime; however the parameter $\varepsilon$ was set to $1$ in this regime, and the same computations show that the residual is actually of size $\mathcal{O}(\varepsilon\mu)$, and therefore of $\mathcal{O}(\mu^2)$ with the Boussinesq scaling regime (\ref{hyp_boussinesqregime}).\quod

\begin{rem}
We remark that for the Serre--Green--Naghdi system, that is without the smallness hypothesis on $\varepsilon$, the situation would be completely different, the expression for the pressure would take a more complex form, incorporating nonlinear effects which would lead to added mass effect for the equation of motion characterizing the solid (for more details we also refer to Section \ref{sec_lastsection}).
\end{rem}

Therefore, following the same computations as in Section \ref{sec_derivationnewton1}, we obtain the same ODE for the solid displacement as (\ref{eq_newtonorder1}), but with a dependence on $\varepsilon$ as well,
\begin{equation}\label{eq_newtonorder2bous}
\ddot{X}_S=-\frac{c_{fric}}{\sqrt{\mu}}\left(\frac{1}{\varepsilon}\tilde{c}_{solid}+\frac{1}{\tilde{M}}\int\limits_{\operatorname{supp}(\mathfrak{b})+ X_S}\zeta\,dx\right)\frac{\dot{X}_S}{\left|\dot{X}_S\right|+\overline{\delta}}+\frac{\varepsilon}{\tilde{M}}\int_{\mathbb{R}^d}\zeta\nabla_x\mathfrak{b}(x-X_S)\,dx.
\end{equation}

Here, we made use of the constant of the solid $\tilde{c}_{solid}$, similar to (\ref{form_solidcoefficient}), defined by
\begin{equation}
\tilde{c}_{solid}=\varepsilon+\frac{|\operatorname{supp}(\mathfrak{b})|}{\tilde{M}}\left(\frac{P_{atm}}{\varrho gH_0}+1\right)-\varepsilon\frac{|\operatorname{Volume}_{Solid}|}{\tilde{M}}.
\end{equation}
The difference between the constants $c_{solid}$ and $\tilde{c}_{solid}$ is the $\varepsilon$ coefficient in the latter one. Due to the additional hypothesis $\varepsilon=1$ in (\ref{hyp_saintvenantregime}), it was not present in the previous section for the Saint-Venant regime, but in the Boussinesq regime (\ref{hyp_boussinesqregime}) it has to be taken into consideration.

Once again, notice the presence of the friction terms in the solid equation, which is potentially of order $(\varepsilon\sqrt{\mu})^{-1}$, so we will have to reason carefully why this doesn't pose a problem for our system. First of all, we have the following concerning the consistency of the solid equation:
\begin{prop}\label{prop_consistorder2}
Let $s_0\geqslant 0$, and let us assume that $\zeta\in\mathcal{C}([0,T];H^{s_0+6}(\mathbb{R}^d))$ and that $\mathfrak{b}\in H^{s_0+6}(\mathbb{R}^d)$ compactly supported. Furthermore let us suppose that $\nabla_x\psi\in \mathcal{C}([0,T];H^{s_0+6}(\mathbb{R}^d))$. In the long wave Boussinesq regime ($\varepsilon = \mathcal{O}(\mu)$) the solid equation (\ref{eq_ndnewton}) is consistent at order $\mathcal{O}(\sqrt{\mu})$ with the model (\ref{eq_newtonorder2bous}) on $[0,T]$ with $T>0$.
\end{prop}
\prove By the regularity assumptions and the additional hypotheses of (\ref{hyp_boussinesqregime}) (Lemmas 3.42. and 5.4. of \cite{lannesbible}), we can write that $\Phi = \psi + \mu^2 R_2$ with
\begin{equation*}
\begin{aligned}
\brck{R_2}_{T,H^{s_0}}&\leqslant C(\brck{\zeta}_{T,H^{s_0+4}}, \|\mathfrak{b}\|_{H^{s_0+4}})\brck{\nabla_x\psi}_{T,H^{s_0+4}}\\
\brck{\partial_tR_2}_{T,H^{s_0}}&\leqslant C(\brck{\zeta}_{T,H^{s_0+6}}, \|\mathfrak{b}\|_{H^{s_0+6}},\brck{\nabla_x\psi}_{T,H^{s_0+6}}),
\end{aligned}
\end{equation*}
meaning that, we have
\begin{equation*}
\begin{aligned}
P|_{z=-1+\varepsilon b}&=\frac{P_{atm}}{\varrho gH_0}+h+\mu^2 R_{P,2}\\
\brck{R_{P,2}}_{T,H^{s_0}}&\leqslant C(\brck{\zeta}_{T,H^{s_0+6}}, \|\mathfrak{b}\|_{H^{s_0+6}},\brck{\nabla_x\psi}_{T,H^{s_0+6}}).
\end{aligned}
\end{equation*}

Hence, in the equation for the solid motion (\ref{eq_ndnewton}), we recover the approximate equation (\ref{eq_newtonorder1}) with the additional error terms
\begin{equation*}
-\sqrt{\mu}\frac{c_{fric}}{\tilde{M}}\frac{\dot{X}_S}{\left|\dot{X}_S\right|+\overline{\delta}}\int_{I(t)}R_{P,2}\,dx+\mu^2\frac{1}{\tilde{M}}\int_{I(t)}R_{P,2}\nabla_x\mathfrak{b}(x-X_S)\,dx,
\end{equation*}
that can be estimated as an $\mathcal{O}(\sqrt{\mu})$ total error term, that is, it is less than
\begin{equation*}
\sqrt{\mu} C(\tilde{M}^{-1},\brck{\zeta}_{T,H^{s_0+6}},\|\mathfrak{b}\|_{H^{s_0+6}},\brck{\nabla_x\psi}_{T,H^{s_0+6}}).
\end{equation*}\quod

We remark that given the fact that the Boussinesq system is consistent at order $\mu^2$ (Corollary 5.20. of \cite{lannesbible}), the consistency of the coupled fluid-solid system can only be at most of order $\sqrt{\mu}$ which is a considerable loss. In order to remedy the situation, we will address some possible extensions of the solid model in Section \ref{sec_lastsection}.

\subsection{The coupled wave-structure model in the Boussinesq regime}\label{sec_boussrewrite}

Here we present some remarks on the right hand side of the Boussinesq system (\ref{eq_boussinesq}). Again, we have that
\begin{equation*}
\pt b(t,x)=-\nabla_x\mathfrak{b}\left(x-X_S(t)\right)\cdot\dot{X}_S(t),
\end{equation*}
however we also have that
\begin{equation*}
\begin{aligned}
\nabla_x\partial_t^2b(t,x)&=\nabla_x\pt\left(-\nabla_x\mathfrak{b}\left(x-X_S(t)\right)\cdot\dot{X}_S(t)\right)\\
&=\nabla_x\left(\nabla_x^2\mathfrak{b}\left(x-X_S(t)\right)\dot{X}_S(t)\cdot\dot{X}_S(t)-\nabla_x\mathfrak{b}\left(x- X_S(t)\right)\cdot\ddot{X}_S(t)\right).
\end{aligned}
\end{equation*}

\textit{To sum it up, the free surface equations with a solid moving at the bottom in the case of the Boussinesq approximation take the following form}
\begin{subequations}
\begin{align}[left = {\,\empheqlbrace}]
\begin{split}
&\partial_t\zeta+\nabla_x\cdot(h\overline{V}) = \pt b,\\
&\left(1-\frac{\mu}{3}\Delta_x\right)\partial_t\overline{V}+\nabla_x\zeta+\varepsilon(\overline{V}\cdot\nabla_x)\overline{V} = -\frac{\mu}{2}\nabla_x\partial_t^2b,
\end{split}  \\
&\ddot{X}_S=-\frac{c_{fric}}{\sqrt{\mu}}\left(\frac{1}{\varepsilon}\tilde{c}_{solid}+\frac{1}{\tilde{M}}\int\limits_{\operatorname{supp}(\mathfrak{b})+ X_S}\zeta\,dx\right)\frac{\dot{X}_S}{\left|\dot{X}_S\right|+\overline{\delta}}+\frac{\varepsilon}{\tilde{M}}\int_{\mathbb{R}^d}\zeta\nabla_x\mathfrak{b}(x-X_S)\,dx.
\end{align}
\end{subequations}

\subsection{A reformulation of the coupled fluid-solid system}

Following the observations of Section \ref{sec_boussrewrite}, we may elaborate the source term of the coupled system. The free surface equations with a solid moving at the bottom in the case of the Boussinesq approximation can be written as
\begin{subequations}\label{eq_systemorder2bous}
\begin{align}[left = {\,\empheqlbrace}]
\begin{split}
&\partial_t\zeta+\nabla_x\cdot(h\overline{V}) = -\nabla_x\mathfrak{b}\left(x-X_S\right)\cdot\dot{X}_S, \\
&\left(1-\frac{\mu}{3}\Delta_x\right)\partial_t\overline{V}+\nabla_x\zeta+\varepsilon(\overline{V}\cdot\nabla_x)\overline{V} = \\
&\qquad-\frac{\mu}{2}\nabla_x\left(\nabla_x^2\mathfrak{b}\left(x-X_S\right)\dot{X}_S\cdot\dot{X}_S-\nabla_x\mathfrak{b}\left(x- X_S\right)\cdot\ddot{X}_S\right),
\end{split}\label{eq_systemorder2bo1} \\
&\ddot{X}_S=-\frac{c_{fric}}{\sqrt{\mu}}\left(\frac{1}{\varepsilon}\tilde{c}_{solid}+\frac{1}{\tilde{M}}\int\limits_{\operatorname{supp}(\mathfrak{b})+X_S}\zeta\,dx\right)\frac{\dot{X}_S}{\left|\dot{X}_S\right|+\overline{\delta}}+\frac{\varepsilon}{\tilde{M}}\int_{\mathbb{R}^d}\zeta\nabla_x\mathfrak{b}(x-X_S)\,dx.\label{eq_systemorder2bo2}
\end{align}
\end{subequations}

First of all, let us remark that a more compact formulation can be derived, just like for the nonlinear Saint-Venant equations coupled with Newton's equation (\ref{eq_systemorder1}) in Section \ref{sec_saintvenantrewrite}. This formula is obtained through the same means as in the previous section, so we will apply similar notations as well. We have the following: the fluid equations (\ref{eq_systemorder2bo1}) for the variable $\mathcal{U}=(\zeta,\overline{V})$ can be written as
\begin{equation}\label{eq_boussinesqsimpl}
D_{\mu}\partial_t\mathcal{U} + \sum_{j=1}^{d}A_j(\mathcal{U},X_S)\partial_j\mathcal{U}+B(\mathcal{U},X_S)=0,
\end{equation}
where the matrix $A_j(\mathcal{U},X_S)$ is the same as the one defined in the previous section, that is
\begin{equation*}
A_j(\mathcal{U},X_S) = 
\left(\begin{array}{c|ccc}
\varepsilon\overline{V}_j & & hI_j & \\
\hline
 & & & \\
I_j^{\top} & & \varepsilon\overline{V}_j\operatorname{Id}_{d\times d} & \\
 & & &
\end{array}\right) \textnormal{ for } 1\leqslant j\leqslant d.
\end{equation*}
We remark that we have the following simple decomposition
\begin{equation}
A_j(\mathcal{U},X_S) = \overline{I}_j + \varepsilon\overline{A}_j(\mathcal{U},X_S) = \left(\begin{array}{c|ccc}
0 & & I_j & \\
\hline
 & & & \\
I_j^{\top} & & 0 & \\
 & & &
\end{array}\right) + \varepsilon\left(\begin{array}{c|ccc}
\overline{V}_j & & (\zeta-b)I_j & \\
\hline
 & & & \\
0 & & \overline{V}_j\operatorname{Id}_{d\times d} & \\
 & & &
\end{array}\right).
\end{equation}

Additionally, we have that
\begin{equation*}
D_{\mu} = 
\left(\begin{array}{c|ccc}
1 & & & \\
\hline
 & & & \\
 & & \left(1-\dfrac{\mu}{3}\Delta_x\right)\operatorname{Id}_{d\times d} & \\
 & & &
\end{array}\right),
\end{equation*}
and the source term vector takes the following form
\begin{equation*}
B(\mathcal{U},X_S)=\begin{pmatrix}
-\varepsilon\overline{V}\cdot\nabla_x\mathfrak{b}\left(x-X_S\right)+\nabla_x\mathfrak{b}\left(x-X_S\right)\cdot\dot{X}_S\\
\dfrac{\mu}{2}\nabla_x\left(\nabla_x^2\mathfrak{b}\left(x-X_S\right)\dot{X}_S\cdot\dot{X}_S-\nabla_x\mathfrak{b}\left(x- X_S\right)\cdot\ddot{X}_S\right)
\end{pmatrix}.
\end{equation*}

\begin{rem}
Once again, we can symmetrize equation (\ref{eq_boussinesqsimpl}) with the use of the matrix
\begin{equation*}
S(\mathcal{U},X_S)=\left(\begin{array}{c|ccc}
1 & & 0 & \\
\hline
 & & & \\
0 & & h\operatorname{Id}_{d\times d} & \\
 & & &
\end{array}\right),
\end{equation*}
remarking that
\begin{equation}
S(\mathcal{U},X_S)=\operatorname{Id}_{(d+1)\times(d+1)}+\varepsilon\overline{S}(\mathcal{U},X_S)=\operatorname{Id}_{(d+1)\times(d+1)}+\varepsilon\left(\begin{array}{c|ccc}
0 & & 0 & \\
\hline
 & & & \\
0 & & (\zeta-b)\operatorname{Id}_{d\times d} & \\
 & & &
\end{array}\right).
\end{equation}
\end{rem}

Let us make one further remark, concerning the second order (nonlinear) ordinary differential equation characterizing the displacement of the solid $X_S$ in (\ref{eq_systemorder2bo2}). Let us adapt the definition of the functional $\mathcal{F}[\mathcal{U}](t,Y,Z)$ introduced in Section \ref{sec_saintvenantrewrite}.
\begin{equation*}
\mathcal{F}[\mathcal{U}](t,Y,Z)=-\frac{c_{fric}}{\sqrt{\mu}}\left(\frac{1}{\varepsilon}\tilde{c}_{solid}+\frac{1}{\tilde{M}}\int\limits_{\operatorname{supp}(b)+Y}\mathcal{U}_0\,dx\right)\frac{Z}{\left|Z\right|+\overline{\delta}}+\frac{\varepsilon}{\tilde{M}}\int_{\mathbb{R}^d}\mathcal{U}_0\nabla_x\mathfrak{b}(x-Y)\,dx.
\end{equation*}

\textit{The coupled system (\ref{eq_systemorder2bous}) has the following equivalent form}
\begin{subequations}\label{eq_systemuniform2}
\begin{align}[left = {\,\empheqlbrace}]
&D_{\mu}\partial_t\mathcal{U} + \sum_{j=1}^{d}A_j(\mathcal{U},X_S)\partial_j\mathcal{U}+B(\mathcal{U},X_S)=0, \label{eq_systemuniform2o1} \\
&\ddot{X}_S = \mathcal{F}[\mathcal{U}]\left(t,X_S,\dot{X}_S\right).\label{eq_systemuniform2o2}
\end{align}
\end{subequations}

\subsection{A priori estimate for the Boussinesq system coupled with Newton's equation}

In this part we present the energy estimate in a Sobolev-type function space for the coupled system (\ref{eq_systemorder2bous}). This estimate is based on classical methods (Gr\"onwall type inequalities), but for an energy functional adapted to the fluid-solid system. In the nonlinear Saint-Venant regime, we constructed an iterative scheme for the system which provided the necessary tools to deduce a local in time existence theorem. The heart of the proof was the energy estimate established on the linearized PDE system (Proposition \ref{prop_structurelinear}) and a separate velocity estimate (Lemma \ref{lemma_velocityestimate}.) for the solid system. Due to the additional dispersive term as well as a more complicated source term on the right hand side of system (\ref{eq_systemorder2bo1}), a refined analysis of the coupling terms is necessary. More precisely the right hand side with $\mu\nabla_x\partial_t^2b$ contains a term of $\mu\ddot{X}_S$ which is asymptotically singular by equation (\ref{eq_systemorder2bo2}).

One additional remark concerns the time of existence of the system. We aim for a long time existence result, which involves the parameter $\varepsilon$. This scale was not present in the previous section since for the Saint-Venant regime (\ref{hyp_saintvenantregime}), we made use of the additional hypothesis of $\varepsilon=1$. However this implies that in the Boussinesq regime (\ref{hyp_boussinesqregime}) more careful estimates are needed; we establish an existence time over a large $\mathcal{O}(\varepsilon^{-1/2})$ scale, while standard methods only provide an $\mathcal{O}(1)$ existence time when the bottom is moving, because of the $\mathcal{O}(1)$ source term $\pt b$ in the first equation of (\ref{eq_boussinesq}). It is however still smaller than the $\mathcal{O}\left(\varepsilon^{-1}\right)$ scale for a fixed bottom (\cite{sautxulong,cosminboussinesq}).

By introducing the wave-structure energy functional
\begin{equation}\label{form_energyvar}
E_B(t)=\frac{1}{2}\int_{\mathbb{R}^d}\zeta^2\,dx+\frac{1}{2}\int_{\mathbb{R}^d}h(\overline{V}\cdot\overline{V})\,dx+\frac{1}{2}\sum_{j=1}^d\int_{\mathbb{R}^d}\frac{\mu}{3}h(\partial_j\overline{V}\cdot\partial_j\overline{V})\,dx+\frac{1}{2\varepsilon}\left|\dot{X}_S\right|^2,
\end{equation}
we can establish first of all an $L^2$ type energy estimate for the coupled system (\ref{eq_systemorder2bous}), from which we will be able to deduce a certain control on the velocity of the solid.

So, we have the following:
\begin{prop}\label{prop_aprioriestimatebous}
Let $\mu\ll 1$ sufficiently small and let us suppose that $s_0>\max(1,d/2)$. Then any $\mathcal{U}\in \mathcal{C}^1([0,T]\times\mathbb{R}^d)\cap\mathcal{C}^1([0,T];H^{s_0}(\mathbb{R}^d))$, $X_S\in\mathcal{C}^2([0,T])$ satisfying the coupled system (\ref{eq_systemuniform2}) (or equivalently (\ref{eq_systemorder2bous})), with initial data $\mathcal{U}(0,\cdot)=\mathcal{U}_{in}\in\mathcal{C}^1(\mathbb{R}^d)\cap H^{s_0}(\mathbb{R}^d)$ and $(X_S(0),\dot{X}_S(0))=(0,v_{S_0})\in\mathbb{R}^d\times\mathbb{R}^d$ verifies the energy estimate
\begin{equation}
\sup_{t\in[0,T]}\left\{e^{-\sqrt{\varepsilon}c_0t}E_B(t)\right\}\leqslant 2E_B(0)+\mu c_0 T\|\mathfrak{b}\|_{H^3}^2,
\end{equation}
where
\begin{equation*}
c_0= c(c_{fric},\tilde{M}^{-1},\brck{\mathcal{U}}_{T,H^{s_0}}, \brck{\mathcal{U}}_{T,W^{1,\infty}},\|\mathfrak{b}\|_{W^{4,\infty}}).
\end{equation*}
\end{prop}

\prove We follow the standard steps of a general energy estimate, adapted for the Boussinesq system with moving bottom, paying close attention to the parameters. We start by multiplying the first equation of (\ref{eq_systemorder2bo1}) by $\zeta$, and the second equation by $h\overline{V}$, after which we integrate on $\mathbb{R}^d$ with respect to the spatial variable $x$. This yields the following system
\begin{align*}[left = {\,\empheqlbrace}]
&\int_{\mathbb{R}^d}\partial_t\zeta\zeta\,dx+\int_{\mathbb{R}^d}\nabla_x\cdot(h\overline{V})\zeta\,dx =  -\int_{\mathbb{R}^d}\zeta\nabla_x\mathfrak{b}\left(x-X_S\right)\,dx\cdot\dot{X}_S, \\
&\int_{\mathbb{R}^d}h\left(1-\frac{\mu}{3}\Delta_x\right)\partial_t\overline{V}\cdot\overline{V}\,dx+\int_{\mathbb{R}^d}h\nabla_x\zeta\overline{V}\,dx+\varepsilon\int_{\mathbb{R}^d}h(\overline{V}\cdot\nabla_x)\overline{V}\cdot\overline{V}\,dx = \\
&\qquad\quad-\frac{\mu}{2}\int_{\mathbb{R}^d}h\nabla_x\left(\nabla_x^2\mathfrak{b}\left(x- X_S\right)\dot{X}_S\cdot\dot{X}_S\right)\cdot\overline{V}\,dx +\frac{\mu}{2}\int_{\mathbb{R}^d}h\nabla_x^2\mathfrak{b}\left(x-X_S\right)\ddot{X}_S\cdot\overline{V}\,dx.
\end{align*}
Our main interest is the terms on the right hand side that represent the coupling in the source term, for the rest we shall reason briefly, since those estimates are part of the classical analysis.

The time derivative term of the second equation can be reformulated by integration by parts in the following way:
\begin{equation*}
\begin{aligned}
\int_{\mathbb{R}^d}h\left(1-\frac{\mu}{3}\Delta_x\right)&\partial_t\overline{V}\cdot\overline{V}\,dx=\frac{1}{2}\frac{d}{dt}\int_{\mathbb{R}^d}h(\overline{V}\cdot\overline{V})\,dx+\frac{1}{2}\frac{d}{dt}\sum_{j=1}^d\int_{\mathbb{R}^d}\frac{\mu}{3}h(\partial_j\overline{V}\cdot\partial_j\overline{V})\,dx \\
&\quad-\frac{1}{2}\int_{\mathbb{R}^d}\pt h\left(\overline{V}\cdot\overline{V}+\frac{\mu}{3}\sum_{j=1}^d(\partial_j\overline{V}\cdot\partial_j\overline{V})\right)\,dx+\frac{\mu}{3}\sum_{j=1}^d\int_{\mathbb{R}^d}\partial_jh (\partial_j\pt\overline{V}\cdot\overline{V})\,dx.
\end{aligned}
\end{equation*}

For the first equation, by making use of an integration by parts as well as equation (\ref{eq_systemorder2bo2}) on the right hand side, we get
\begin{align*}[left = {\,\empheqlbrace}]
&\frac{1}{2}\frac{d}{dt}\int_{\mathbb{R}^d}\zeta^2\,dx+\frac{\tilde{M}}{2}\frac{1}{\varepsilon}\frac{d}{dt}\left|\dot{X}_S\right|^2+\int_{\mathbb{R}^d}\varepsilon\nabla_x\zeta\cdot\overline{V}\zeta\,dx+\int_{\mathbb{R}^d}h(\nabla_x\cdot\overline{V})\zeta\,dx =\\
&\qquad\qquad \varepsilon\int_{\mathbb{R}^d}\nabla_x\mathfrak{b}\left(x-X_S\right)\cdot\overline{V}\zeta\,dx-\frac{\tilde{M} c_{fric}}{\varepsilon\sqrt{\mu}}\left(\frac{1}{\varepsilon}\tilde{c}_{solid}+\frac{1}{\tilde{M}}\int_{I(t)}\zeta\,dx\right)\frac{\left|\dot{X}_S\right|^2}{\left|\dot{X}_S\right|+\overline{\delta}}, \\
&\frac{1}{2}\frac{d}{dt}\int_{\mathbb{R}^d}h(\overline{V}\cdot\overline{V})\,dx+\frac{1}{2}\frac{d}{dt}\sum_{j=1}^d\int_{\mathbb{R}^d}\frac{\mu}{3}h(\partial_j\overline{V}\cdot\partial_j\overline{V})\,dx \\
&\qquad-\frac{1}{2}\int_{\mathbb{R}^d}\pt h\left(\overline{V}\cdot\overline{V}+\frac{\mu}{3}\sum_{j=1}^d(\partial_j\overline{V}\cdot\partial_j\overline{V})\right)\,dx+\frac{\mu}{3}\sum_{j=1}^d\int_{\mathbb{R}^d}\partial_jh (\partial_j\pt\overline{V}\cdot\overline{V})\,dx\\
&\qquad+\int_{\mathbb{R}^d}h\nabla_x\zeta\overline{V}\,dx+\varepsilon\int_{\mathbb{R}^d}h(\overline{V}\cdot\nabla_x)\overline{V}\cdot\overline{V}\,dx = \\
&\qquad\qquad-\frac{\mu}{2}\int_{\mathbb{R}^d}h\nabla_x\left(\nabla_x^2\mathfrak{b}\left(x- X_S\right)\dot{X}_S\cdot\dot{X}_S\right)\cdot\overline{V}\,dx +\frac{\mu}{2}\int_{\mathbb{R}^d}h\nabla_x^2\mathfrak{b}\left(x-X_S\right)\ddot{X}_S\cdot\overline{V}\,dx.
\end{align*}
Notice that by equation (\ref{eq_systemorder2bo2}), we have been able to substitute part of the contribution associated to the source term $\pt b$ as a component of the energy $E_B(t)$ on the left hand side of the first equation. This is crucial to get an extended existence time. Moreover, on the right hand side, a now nonpositive friction term appeared that can be easily controlled.

Now we add together these two equations and in what follows, by making use of term by term estimates, we arrive to a Gr\"onwall-type inequality concerning the energy functional $E_B(t)$ (for $0\leqslant t\leqslant T$) which then allows us to properly conclude the demonstration. Hence we are left with
\begin{equation}
\frac{d}{dt}E_B(t) = A_B + B_B + C_B + D_B + F_B + G_B,
\end{equation}
where
\begin{align}
A_B &:= \frac{1}{2}\int_{\mathbb{R}^d}\pt h\left(\overline{V}\cdot\overline{V}+\frac{\mu}{3}\sum_{j=1}^d(\partial_j\overline{V}\cdot\partial_j\overline{V})\right)\,dx,\\
B_B &:= -\frac{\mu}{3}\sum_{j=1}^d\int_{\mathbb{R}^d}\partial_jh (\partial_j\pt\overline{V}\cdot\overline{V})\,dx,\\
C_B &:= -\int_{\mathbb{R}^d}\varepsilon\nabla_x\zeta\cdot\overline{V}\zeta\,dx-\int_{\mathbb{R}^d}h(\nabla_x\cdot\overline{V})\zeta\,dx-\int_{\mathbb{R}^d}h\nabla_x\zeta\overline{V}\,dx-\varepsilon\int_{\mathbb{R}^d}h(\overline{V}\cdot\nabla_x)\overline{V}\cdot\overline{V}\,dx,\\
D_B &:= \varepsilon\int_{\mathbb{R}^d}\nabla_x\mathfrak{b}\left(x-X_S\right)\cdot\overline{V}\zeta\,dx-\frac{\tilde{M} c_{fric}}{\varepsilon\sqrt{\mu}}\left(\frac{1}{\varepsilon}\tilde{c}_{solid}+\frac{1}{\tilde{M}}\int_{I(t)}\zeta\,dx\right)\frac{\left|\dot{X}_S\right|^2}{\left|\dot{X}_S\right|+\overline{\delta}}, \\
F_B &:= -\frac{\mu}{2}\int_{\mathbb{R}^d}h\nabla_x\left(\nabla_x^2\mathfrak{b}\left(x- X_S\right)\dot{X}_S\cdot\dot{X}_S\right)\cdot\overline{V}\,dx, \\
G_B &:= \frac{\mu}{2}\int_{\mathbb{R}^d}h\nabla_x^2\mathfrak{b}\left(x-X_S\right)\ddot{X}_S\cdot\overline{V}\,dx.
\end{align}
 
Now we proceed to estimate each term on the right hand side. By making use of the first equation of the Boussinesq system (\ref{eq_systemorder2bous}), namely that
\begin{equation*}
\pt h = \varepsilon\nabla_x\cdot(h\overline{V}),
\end{equation*}
we can establish that
\begin{equation*}
A_B\leqslant\varepsilon c(\|\mathcal{U}\|_{W^{1,\infty}},\|\mathfrak{b}\|_{W^{1,\infty}})\left(\|\overline{V}\|_{L^2}^2+\frac{\mu}{3}\|\overline{V}\|_{H^1}^2\right).
\end{equation*}

As for the term $B_B$, we aim to estimate the $L^2$ norm of the mixed derivative term $\partial_j\pt\overline{V}$. By making use of the second equation of the system (\ref{eq_systemorder2bous}), we have
\begin{align}\label{id_termmixte}
\frac{\mu}{3}\partial_j\pt\overline{V}&=\left(1-\frac{\mu}{3}\Delta_x\right)^{-1}\left(-\frac{\mu}{3}\partial_j\nabla_x\right)\zeta+\left(1-\frac{\mu}{3}\Delta_x\right)^{-1}\varepsilon\left(-\frac{\mu}{3}\partial_j((\overline{V}\cdot\nabla_x)\overline{V})\right)\\
&\qquad+\frac{\mu}{2}\left(1-\frac{\mu}{3}\Delta_x\right)^{-1}\left(-\frac{\mu}{3}\partial_j\nabla_x\right)\pt^2b.\nonumber
\end{align}
Let us estimate each term separately. Given the fact that $\left(1-\frac{\mu}{3}\Delta_x\right)^{-1}\left(-\frac{\mu}{3}\partial_j\nabla_x\right)$ is a zeroth order differential operator whose symbol is uniformly bounded with respect to $\mu$, we can easily deduce that
\begin{equation*}
\left\|\left(1-\frac{\mu}{3}\Delta_x\right)^{-1}\left(-\frac{\mu}{3}\partial_j\nabla_x\right)\zeta\right\|_{L^2}\lesssim\|\zeta\|_{L^2}.
\end{equation*}
For the second term, first of all, we have that the operator $\left(1-\frac{\mu}{3}\Delta_x\right)^{-1}\mu\partial_j$ has a symbol of order $-1$, uniformly bounded with respect to $\mu$, therefore
\begin{equation*}
\left\|\left(1-\frac{\mu}{3}\Delta_x\right)^{-1}\varepsilon\left(-\frac{\mu}{3}\partial_j((\overline{V}\cdot\nabla_x)\overline{V})\right)\right\|_{L^2}\lesssim\varepsilon\|(\overline{V}\cdot\nabla_x)\overline{V})\|_{H^{-1}},
\end{equation*}
from which, by a classical product estimate, we have that for $s_0\geqslant d/2$, $-1>-s_0$
\begin{equation*}
\|(\overline{V}\cdot\nabla_x)\overline{V})\|_{H^{-1}}\lesssim\|\overline{V}\|_{H^{s_0}}\|\nabla_x\cdot\overline{V}\|_{H^{-1}}\lesssim\|\overline{V}\|_{H^{s_0}}\|\overline{V}\|_{L^2}.
\end{equation*}
As for the third term from (\ref{id_termmixte}), we use the chain rule for $\pt^2b$, as well as the fact that $\left(1-\frac{\mu}{3}\Delta_x\right)^{-1}$ is uniformly bounded in $\mu$ as an differential operator of order $0$. This yields
\begin{equation*}
\mu^2\left\|\left(1-\frac{\mu}{3}\Delta_x\right)^{-1}\left(-\partial_j\nabla_x\right)\pt^2b\right\|_{L^2}\lesssim \mu^2(\|\mathfrak{b}\|_{W^{4,\infty}}|\dot{X}_S|^2+\|\mathfrak{b}\|_{H^3}|\ddot{X}_S|),
\end{equation*}
where we made use of the fact that $\mathfrak{b}$ is compactly supported.  Here we can estimate $\mu^2 |\ddot{X}_S|$ directly from equation (\ref{eq_systemorder2bo2}) due to the additional smallness parameter. More exactly we have that
\begin{equation*}
\mu^2\frac{c_{fric}}{\sqrt{\mu}\varepsilon}\tilde{c}_{solid}\leqslant\sqrt{\mu}C(c_{fric},\|\mathfrak{b}\|_{L^\infty})\quad\textnormal{ and }\quad\frac{|\dot{X}_S|}{|\dot{X}_S|+\overline{\delta}}\leqslant 1,
\end{equation*}
which allows us to infer that
\begin{equation*}
\mu^2|\ddot{X}_S|\leqslant\sqrt{\mu}C(c_{fric},\|\mathfrak{b}\|_{W^{1,\infty}},\|\zeta\|_{L^\infty})
\end{equation*}

To sum it up, we have obtained the following estimate:
\begin{equation}\label{ineq_gradvelobouss}
\left\|\frac{\mu}{3}\partial_j\pt\overline{V}\right\|_{L^2}\lesssim\|\zeta\|_{L^2}+\varepsilon c(\|\mathcal{U}\|_{H^{s_0}})\|\overline{V}\|_{L^2}+\mu^2|\dot{X}_S|^2+\sqrt{\mu}\|\mathfrak{b}\|_{H^3}.
\end{equation}
Thus we get
\begin{equation*}
B_B\leqslant \varepsilon c(\|\mathcal{U}\|_{W^{1,\infty}},\|\mathfrak{b}\|_{W^{4,\infty}})\left[\|\zeta\|_{L^2}\|\overline{V}\|_{L^2}+\varepsilon c(\|\mathcal{U}\|_{H^{s_0}})\|\overline{V}\|_{L^2}^2+\left|\dot{X}_S\right|^2+\sqrt{\mu}\|\mathfrak{b}\|_{H^3}\|\overline{V}\|_{L^2}\right],
\end{equation*}
here the last term can be estimated as $\|\nabla_x\mathfrak{b}\|_{L^2}^2+\|\overline{V}\|_{L^2}^2$ as well.

The integrals incorporating the nonlinear spatial derivative terms correspond to
\begin{equation*}
C_B = \sum_{j=1}^{d}\int_{\mathbb{R}^d}SA_j(\mathcal{U})\partial_j\mathcal{U}\cdot\mathcal{U}\,dx=-\frac{1}{2}\sum_{j=1}^{d}\int_{\mathbb{R}^d}\partial_j(SA_j)(\mathcal{U})\mathcal{U}\cdot\mathcal{U}\,dx,
\end{equation*}
to which we can easily find an upper bound, giving
\begin{equation*}
C_B\leqslant \varepsilon c(\|\mathcal{U}\|_{W^{1,\infty}},\|\mathfrak{b}\|_{W^{1,\infty}})\|\mathcal{U}\|_{L^2}^2.
\end{equation*}

For the first two source terms for the system, basic $L^\infty$-norm estimates and Cauchy-Schwartz inequalities provide the necessary means to conclude
\begin{align*}
D_B&\leqslant\varepsilon c(\|\mathfrak{b}\|_{W^{1,\infty}})\|\zeta\|_{L^2}\|\overline{V}\|_{L^2} + 0,\\
F_B&\leqslant\mu c(\|\mathcal{U}\|_{W^{1,\infty}},\|\mathfrak{b}\|_{W^{3,\infty}})\left|\dot{X}_S\right|^2.
\end{align*}
We remark that the friction term can be straightforwardly bounded above by $0$ in the estimate for $D_B$.

We leave the last source term, $G_B$, as it is due to the presence of $\ddot{X}_S(t)$; according to equation (\ref{eq_systemorder2bo2}), it requires some attention to avoid problems arising from the friction part (the asymptotically singular terms).

So, to sum up the previous estimates, we get that
\begin{align*}
\frac{d}{dt}E_B(t)\leqslant &+\varepsilon c(\|\mathcal{U}\|_{W^{1,\infty}},\|\mathfrak{b}\|_{W^{1,\infty}})\left(\|\overline{V}\|_{L^2}^2+\frac{\mu}{3}\|\nabla_x\cdot\overline{V}\|_{L^2}^2\right) \\
&+\varepsilon c(\|\mathcal{U}\|_{W^{1,\infty}},\|\mathfrak{b}\|_{W^{1,\infty}})\left[\|\zeta\|_{L^2}\|\overline{V}\|_{L^2}+\varepsilon c(\|\mathcal{U}\|_{H^{s_0}})\|\overline{V}\|_{L^2}^2 +\sqrt{\mu}\|\overline{V}\|_{L^2}^2 +\|\mathcal{U}\|_{L^2}^2\right]\\
&+\varepsilon c(\|\mathcal{U}\|_{W^{1,\infty}},\|\mathfrak{b}\|_{W^{4,\infty}})\|\mathcal{U}\|_{L^2}^2(\mu^2|\dot{X}_S|^2+\sqrt{\mu}\|\mathfrak{b}\|_{H^3}^2) \\
&+\varepsilon c(\|\mathfrak{b}\|_{W^{1,\infty}})\|\zeta\|_{L^2}\|\overline{V}\|_{L^2} + 0 + \mu c(\|\mathcal{U}\|_{W^{1,\infty}},\|\mathfrak{b}\|_{W^{3,\infty}})\left|\dot{X}_S\right|^2\\ 
&+ G_B.
\end{align*}

So we may deduce that
\begin{equation}
\frac{d}{dt}E_B(t)\leqslant \varepsilon c(\|\mathcal{U}\|_{H^{s_0}}, \|\mathcal{U}\|_{W^{1,\infty}},\| \mathfrak{b}\|_{W^{4,\infty}})\left(E_B(t) +\sqrt{\mu} \|\mathfrak{b}\|_{H^3}^2\right) + G_B
\end{equation}
from which by integrating with respect to the time variable $t$ (keeping in mind that $0\leqslant t\leqslant T$), we obtain
\begin{equation}\label{ineq_pregronwall}
E_B(t)-E_B(0)\leqslant \varepsilon c_0\int_0^tE_B(\tau)\,d\tau +\sqrt{\mu}\varepsilon c_0t\|\mathfrak{b}\|_{H^3}^2+\int_0^t G_B\,d\tau,
\end{equation}
where we made use of the constant
\begin{equation*}
c_0= c(\brck{\mathcal{U}}_{T,H^{s_0}}, \brck{\mathcal{U}}_{T,W^{1,\infty}},\|\mathfrak{b}\|_{W^{4,\infty}}).
\end{equation*}

\begin{lem}\label{lem_derniersourceterm}
The remaining source term $G_B$ satisfies the following estimate for all $0\leqslant t\leqslant T$,
\begin{equation}
\begin{aligned}
\frac{\mu}{2}\int_0^t \int_{\mathbb{R}^d}h\nabla_x^2\mathfrak{b}\left(x-X_S\right)\ddot{X}_S\cdot\overline{V}\,dx\,d\tau\leqslant&\mu c(\|\mathfrak{b}\|_{W^{2,\infty}}, \brck{\mathcal{U}}_{T,L^\infty})(E_B(t)+E_B(0))\\
&+\sqrt{\varepsilon}c(c_{fric},\tilde{M}^{-1},\|\mathfrak{b}\|_{W^{4,\infty}}, \brck{\mathcal{U}}_{T,H^{s_0}}, \brck{\mathcal{U}}_{T,W^{1,\infty}})\int_0^tE_B(\tau)\,d\tau \\
&+\mu c(\|\mathfrak{b}\|_{W^{4,\infty}}, \brck{\mathcal{U}}_{T,W^{1,\infty}})\int_0^t \|\mathfrak{b}\|_{H^3}^2\,d\tau.
\end{aligned}
\end{equation}
\end{lem}

\prove To handle the source term, first of all we apply an integration by parts in the time variable. This yields
\begin{equation*}
\begin{aligned}
\mu \int_{\mathbb{R}^d}\int_0^t &h\nabla_x^2\mathfrak{b}\left(x- X_S\right)\ddot{X}_S\cdot\overline{V}\,d\tau\,dx =\mu \int_{\mathbb{R}^d}h\nabla_x^2\mathfrak{b}\left(x-X_S\right)\dot{X}_S\cdot\overline{V}(t,x)\,dx\\ 
&-\mu \int_{\mathbb{R}^d}h\nabla_x^2\mathfrak{b}\left(x-0\right)v_{S_0}\cdot\overline{V}_{in}(x)\,dx-\mu \int_{\mathbb{R}^d}\partial_th\nabla_x^2\mathfrak{b}\left(x-X_S\right)\dot{X}_S\cdot\overline{V}(t,x)\,dx\\
&+\mu \int_0^t\int_{\mathbb{R}^d}h\nabla_x(\nabla_x^2\mathfrak{b}\left(x- X_S(\tau)\right)\dot{X}_S(\tau)\cdot\dot{X}_S(\tau))\cdot\overline{V}\,dx\,d\tau\\
&-\mu \int_0^t\int_{\mathbb{R}^d}h\nabla_x^2\mathfrak{b}\left(x- X_S(\tau)\right)\dot{X}_S(\tau)\cdot\partial_t\overline{V}\,dx\,d\tau.
\end{aligned}
\end{equation*}
The first two boundary terms can be estimated similarly, for the first term we have that
\begin{equation*}
\mu \int_{\mathbb{R}^d}h\nabla_x^2\mathfrak{b}\left(x- X_S\right)\dot{X}_S\cdot\overline{V}(t,x)\,dx\leqslant\mu c(\|\mathfrak{b}\|_{W^{2,\infty}},\|\mathcal{U}\|_{L^\infty})(|\dot{X}_S|^2+\|\overline{V}\|_{L^2}^2),
\end{equation*}
and we can deduce an identical estimate for the initial data. Here we used the fact that $\mathfrak{b}$ is compactly supported, and as such the integrals can be calculated on $\operatorname{supp}(\mathfrak{b})$. Since we assume that $\mu$ is sufficiently small, this estimate with the energy term will be absorbed by the energy term on the left hand side of (\ref{ineq_pregronwall}).

Once again making use of the first equation of (\ref{eq_systemorder2bous}) we obtain
\begin{equation*}
\mu \int_{\mathbb{R}^d}\partial_th\nabla_x^2\mathfrak{b}\left(x-X_S\right)\dot{X}_S\cdot\overline{V}(t,x)\,dx\leqslant \mu\varepsilon c(\|\mathfrak{b}\|_{W^{2,\infty}},\|\mathcal{U}\|_{W^{1,\infty}}) (|\dot{X}_S|^2+\|\overline{V}\|_{L^2}^2).
\end{equation*}

The integral on the support of $\mathfrak{b}$ gives
\begin{equation*}
\begin{aligned}
\mu \int_{\mathbb{R}^d}h\nabla_x(\nabla_x^2\mathfrak{b}\left(x- X_S(\tau)\right)\dot{X}_S\cdot\dot{X}_S)\cdot\overline{V}\,dx\,d\tau&\leqslant\mu c(\|\mathfrak{b}\|_{L^{\infty}},\|\mathcal{U}\|_{L^\infty})|X_S|^2\|\nabla_x^3\mathfrak{b}\|_{L^2(\operatorname{supp}(\mathfrak{b}))}\|\overline{V}\|_{L^2(\operatorname{supp}(\mathfrak{b}))} \\
&\leqslant \mu c(\|\mathfrak{b}\|_{W^{3,\infty}},\|\mathcal{U}\|_{L^\infty})|X_S|^2.
\end{aligned}
\end{equation*}

Finally, by an integration by parts with respect to the spatial variable, we get that
\begin{equation*}
\int_0^t\int_{\mathbb{R}^d}\nabla_x^2\mathfrak{b}\left(x- X_S(\tau)\right)\dot{X}_S(\tau)\cdot\partial_t\overline{V}\,dx\,d\tau= 
\int_0^t\int_{\mathbb{R}^d}\nabla_x\mathfrak{b}\left(x- X_S(\tau)\right)\cdot\dot{X}_S(\tau)(\nabla_x\cdot\partial_t\overline{V})\,dx\,d\tau,
\end{equation*}
from which we deduce that
\begin{equation*}
\begin{aligned}
\mu \int_0^t\int_{\mathbb{R}^d}&h\nabla_x^2\mathfrak{b}\left(x- X_S(\tau)\right)\dot{X}_S\cdot\partial_t\overline{V}\,dx\,d\tau\leqslant  c(\|\mathfrak{b}\|_{W^{2,\infty}}, \brck{\mathcal{U}}_{T,L^\infty})\int_0^t|\dot{X}_S|\left\|\dfrac{\mu}{3}\nabla_x\cdot\partial_t\overline{V}\right\|_{L^2}\,dx\,d\tau\\
&\leqslant\sqrt{\varepsilon} c(\|\mathfrak{b}\|_{W^{2,\infty}}, \brck{\mathcal{U}}_{T,L^\infty})\int_0^t\frac{|\dot{X}_S(\tau)|}{\sqrt{\varepsilon}}\left(\|\zeta\|_{L^2}+\varepsilon c(\|\mathcal{U}\|_{H^{s_0}})\|\overline{V}\|_{L^2}+\sqrt{\mu}\|\overline{V}\|_{L^2}+\mu|\dot{X}_S|\right)\,d\tau,\\
&\;+\mu c(\|\mathfrak{b}\|_{W^{1,\infty}}, \brck{\mathcal{U}}_{T,W^{1,\infty}})\int_0^t \|\mathfrak{b}\|_{H^3}^2\,d\tau
\end{aligned}
\end{equation*}
where we made use of our previous observation adapted to $\mu\nabla_x\cdot\partial_t\overline{V}$ (inequality (\ref{ineq_gradvelobouss})). Remarking that $\varepsilon^{-1/2}|X_S|\leqslant E_B^{1/2}$ by definition, we obtain
\begin{equation*}
\begin{aligned}
\mu \int_0^t\int_{\mathbb{R}^d}&h\nabla_x^2\mathfrak{b}\left(x- X_S(\tau)\right)\dot{X}_S(\tau)\cdot\partial_t\overline{V}\,dx\,d\tau\leqslant\\
&\sqrt{\varepsilon}c(c_{fric},\tilde{M}^{-1}, \|\mathfrak{b}\|_{W^{2,\infty}}, \brck{\mathcal{U}}_{T,H^{s_0}})\int_0^tE_B(\tau)\,d\tau+\mu c(\|\mathfrak{b}\|_{W^{1,\infty}},\brck{\mathcal{U}}_{T,W^{1,\infty}})\int_0^t \|\mathfrak{b}\|_{H^3}^2\,d\tau,
\end{aligned}
\end{equation*}
which in turn allows us to conclude this lemma.\quod

So by Lemma \ref{lem_derniersourceterm}. and inequality (\ref{ineq_pregronwall}), we obtain that for $\mu$ sufficiently small
\begin{equation}
E_B(t)\leqslant 2E_B(0)+\sqrt{\varepsilon}\tilde{c}_0\int_0^tE_B(\tau)\,d\tau +\mu \tilde{c}_0t \|\mathfrak{b}\|_{H^3}^2,
\end{equation}
with the constant
\begin{equation*}
\tilde{c}_0= c(c_{fric},\tilde{M}^{-1},\brck{\mathcal{U}}_{T,H^{s_0}}, \brck{\mathcal{U}}_{T,W^{1,\infty}},\|\mathfrak{b}\|_{W^{4,\infty}}).
\end{equation*}

Thus, by Gr\"onwall's inequality, we can conclude the energy estimate. \quod

This concludes the $L^2$-estimates (case $s=0$). Let us mention some consequences concerning the velocity of the solid.

\begin{cor}\label{cor_velocityestimate} This energy estimate provides us with a natural control on the solid velocity, namely
\begin{equation}\label{form_veloestim}
\sup_{t\in[0,T]}\left\{e^{-\sqrt{\varepsilon} c_0t}\left|\dot{X}_S(t)\right|^2\right\}\leqslant\varepsilon\|\mathcal{U}_{in}\|_{\mathcal{X}^0}^2+\left|v_{S_0}\right|^2+\varepsilon\mu c_0T \|\mathfrak{b}\|_{H^3}^2,
\end{equation}
where
\begin{equation*}
c_0= c(c_{fric},\tilde{M}^{-1},\brck{\mathcal{U}}_{T,W^{1,\infty}},\brck{\mathcal{U}}_{T,H^{s_0}},\|\mathfrak{b}\|_{W^{4,\infty}}).
\end{equation*}
This implies that the solid velocity stays bounded on a $\mathcal{O}(\varepsilon^{-1/2})$ timescale as long as $c_0$ stays bounded.
\end{cor}

\begin{rem}
Following the steps of Lemma \ref{lemma_velocityestimate}. we would have obtained the velocity estimate
\begin{equation}
\left|\dot{X}(t)\right|\leqslant \left|v_{S_0}\right|+\varepsilon\frac{\|\mathcal{U}\|_{L^\infty}\|b\|_{W^{1,\infty}}}{\tilde{M}}t,
\end{equation}
which is a worse estimate than the one presented in the previous corollary and it cannot be used to obtain an extended existence time.
\end{rem}

\begin{rem}\label{rem_existencetimesolid}
By the identity
\begin{equation*}
\dot{X}_S(t) = \sqrt{\varepsilon}\frac{\dot{X}_S(t)}{\sqrt{\varepsilon}},
\end{equation*}
from (\ref{form_veloestim}) of Corollary \ref{cor_velocityestimate}, it is easy to see that if the initial velocity is of order $\mathcal{\sqrt{\varepsilon}}$, that is $\varepsilon^{-1/2}v_{S_0}$ is uniformly bounded in $\mu$ and $\varepsilon$, then the scaled solid velocity $\varepsilon^{-1/2}\dot{X}_S(t)$ stays uniformly bounded. Moreover, this uniform bound is valid up until a time of order $\mathcal{O}(\varepsilon^{-1/2})$ as long as $c_0$ remains bounded.
\end{rem}

For higher order energy estimates we are going to make use of this estimate and the differential operator $\Lambda^s=(1-\Delta_x)^{s/2}$. The energy functional associated to these estimates writes as
\begin{equation*}
E_B^s(t)=\frac{1}{2}\int_{\mathbb{R}^d}(\Lambda^s\zeta)^2\,dx+\frac{1}{2}\int_{\mathbb{R}^d}h(\Lambda^s\overline{V}\cdot\Lambda^s\overline{V})\,dx+\frac{1}{2}\sum_{j=1}^d\int_{\mathbb{R}^d}\frac{\mu}{3}h(\partial_j\Lambda^s\overline{V}\cdot\partial_j\Lambda^s\overline{V})\,dx+\frac{1}{2\varepsilon}\left|\dot{X}_S\right|^2.
\end{equation*}

Due to the special structure of our system, let us define the following adapted Sobolev space to provide a uniformly formulated energy estimate.
\begin{defi}
The Sobolev-type space $\mathcal{X}^s$ is given by
\begin{equation*}
\mathcal{X}^s(\mathbb{R}^d)=\left\{\mathcal{U}=(\zeta,\overline{V})\in L^2(\mathbb{R}^d)\textnormal{ such that }\|\mathcal{U}\|_{\mathcal{X}^s}<\infty\right\},
\end{equation*}
where
\begin{equation*}
\|\mathcal{U}\|_{\mathcal{X}^s}=\|\zeta\|_{H^s}+\|\overline{V}\|_{H^s}+\sqrt{\mu}\|\overline{V}\|_{H^{s+1}}.
\end{equation*}
\end{defi}
The last term in the $\mathcal{X}^s$ norm appeared due to the necessity to control the dispersive smoothing through $\sqrt{\mu}$ times the partial derivatives.

We have to modify certain parts of the proof, due to the fact that some of the cancellations used above cease to work anymore. More precisely, we have that

\begin{prop}\label{prop_aprioriestimatebous2}
Let $\mu\ll 1$ sufficiently small and let us take $s\in\mathbb{R}$ with $s>d/2+1$. Let us take $\mathcal{U}\in \mathcal{C}([0,T];\mathcal{X}^s(\mathbb{R}^d))\cap C^1([0,T];\mathcal{X}^{s-1}(\mathbb{R}^d))$, $X_S\in\mathcal{C}^1([0,T])$ satisfying the coupled system (\ref{eq_systemuniform2}) (or equivalently (\ref{eq_systemorder2bous})), with initial data $\mathcal{U}(0,\cdot)\in\mathcal{X}^s(\mathbb{R}^d)$ and $(X_S(0),\dot{X}_S(0))=(0,\sqrt{\varepsilon}V_{S_0})\in\mathbb{R}^d\times\mathbb{R}^d$. Then $\mathcal{U}$, $X_S$ verifies the energy estimate
\begin{equation}
\sup_{t\in[0,T]}\left\{e^{-\sqrt{\varepsilon} c_{s}t}\left(\frac{1}{2}\|\mathcal{U}\|_{\mathcal{X}^s}^2+\frac{1}{2\varepsilon}\left|\dot{X}_S\right|^2\right)\right\}\leqslant 2\|\mathcal{U}(0,\cdot)\|_{\mathcal{X}^s}^2+2\left|V_{S_0}\right|^2+\sqrt{\varepsilon}Tc_s\|\mathfrak{b}\|_{H^{s+3}}^2,
\end{equation}
where
\begin{equation*}
c_{s}= c(c_{fric},\tilde{M}^{-1},\brck{\mathcal{U}}_{T,H^s},\|\mathfrak{b}\|_{H^{s+3}}).
\end{equation*}
\end{prop}

\begin{rem}\label{rem_existencetimevelo}
Notice that taking into account the coupling effect for the a priori estimate ensured that the constant in the exponential stays of order $\sqrt{\varepsilon}$, which guarantees a proper control on the fluid velocity over a time $\mathcal{O}(\varepsilon^{-1/2})$, which is better than what the general theory would imply for a time dependent bottom variation.
\end{rem}

\prove We start by applying the operator $\Lambda^s$ on the symmetrized equation ((\ref{eq_systemuniform2}) multiplied by $S(\mathcal{U},X_S)$), and we would like to use the techniques presented for the case of $s=0$, treating $\Lambda^s\mathcal{U}$ as our new unknown. Thus we are left with
\begin{equation}\label{eq_boussinesqcommutators}
\begin{aligned}
S(\mathcal{U},X_S)D_{\mu}\partial_t\Lambda^s\mathcal{U} &+ \sum_{j=1}^{d}SA_j(\mathcal{U},X_S)\partial_j\Lambda^s\mathcal{U}+\Lambda^sSB(\mathcal{U},X_S)\\
&+[\Lambda^s,S(\mathcal{U},X_S)]D_{\mu}\partial_t\mathcal{U}+\sum_{j=1}^{d}[\Lambda^s,SA_j(\mathcal{U},X_S)]\partial_j\mathcal{U}=0.
\end{aligned}
\end{equation}
Notice the presence of the additional commutator terms in the equation.

Our main idea is the same as before, after multiplying the equation by $\Lambda^s\mathcal{U}$ and integrating over $\mathbb{R}^d$, we make use of similar estimates as in the first part for the $L^2$ estimate to obtain a Gr\"onwall type inequality for the corresponding modified energy functional $E_B^s(t)$.

For the first two terms of our new equation, which correspond to the time derivative and nonlinear terms of the original equation (\ref{eq_systemsymmorder2reg1}), they may be treated similarly as before, obtaining the same estimates with the same constants, only for $H^s$-norm instead of $L^2$-norm. 

The main difference is the presence of the commutators in equation (\ref{eq_boussinesqcommutators}) due to $\Lambda^s$, and the treatment of the source term since the cancellation obtained by using the ODE (\ref{eq_systemorder2bo2}) does not work anymore. We will make use of the well-known Kato-Ponce inequality (for this we have $s>0$) as well as Sobolev-embedding results (for these, the condition $s>d/2+1$ is necessary) to establish commutator estimates. Namely, we have that for $f\in H^s$, $g\in H^{s-1}$
\begin{equation*}
\|[\Lambda_s,f]g\|_{L^2}\lesssim\|f\|_{H^s}\|g\|_{L^\infty}+\|f\|_{W^{1,\infty}}\|g\|_{H^{s-1}}\lesssim\|f\|_{H^s}\|g\|_{H^{s-1}},
\end{equation*}
the latter inequality coming from the embedding $H^s(\mathbb{R}^d)\hookrightarrow W^{1,\infty}(\mathbb{R}^d)$. We also have that for $f\in H^s$, $g\in H^s$, and $s>d/2$
\begin{equation*}
\|\Lambda_s(fg)\|_{L^2}\lesssim\|f\|_{H^s}\|g\|_{H^s}.
\end{equation*}

So, by the decomposition of the symmetrizer matrix, we may write that
\begin{align*}
&\left|\int_{\mathbb{R}^d}[\Lambda^s,S(\mathcal{U},X_S)]D_{\mu}\partial_t\mathcal{U}\cdot\Lambda^s\mathcal{U}\,dx\right|=\left|\int_{\mathbb{R}^d}[\Lambda^s,\varepsilon\overline{S}(\mathcal{U},X_S)]D_{\mu}\partial_t\mathcal{U}\cdot\Lambda^s\mathcal{U}\,dx\right|\\
&\quad\lesssim\varepsilon(\|\mathcal{U}\|_{H^s}+\|\mathfrak{b}\|_{H^s})\left\|\sum_{j=1}^{d}SA_j(\mathcal{U},X_S)\partial_j\mathcal{U}+SB(\mathcal{U},X_S)\right\|_{H^{s-1}}\cdot\|\mathcal{U}\|_{H^s}\\
&\quad\lesssim \varepsilon c(\|\mathcal{U}\|_{H^s},\|\mathfrak{b}\|_{H^s})\|\mathcal{U}\|_{H^s}^2+\varepsilon c(\|\mathcal{U}\|_{H^s},\|\mathfrak{b}\|_{H^s})\|SB(\mathcal{U},X_S)\|_{H^{s}}\|\mathcal{U}\|_{H^s}.
\end{align*}
Here we made use of equation (\ref{eq_systemuniform2}) to handle the time derivative. The additional term will be absorbed by the source term in equation (\ref{eq_boussinesqcommutators}).

Additionally we have that
\begin{align*}
&\left|\int_{\mathbb{R}^d}\sum_{j=1}^{d}[\Lambda^s,SA_j(\mathcal{U},X_S)]\partial_j\mathcal{U}\cdot\Lambda^s\mathcal{U}\,dx\right|=\left|\int_{\mathbb{R}^d}\sum_{j=1}^{d}[\Lambda^s,\varepsilon S\overline{A}_j(\mathcal{U},X_S)]\partial_j\mathcal{U}\cdot\Lambda^s\mathcal{U}\,dx\right|\\
&\quad\lesssim\varepsilon\sum_{j=1}^{d}\left\|S\overline{A}_j(\mathcal{U},X_S)\right\|_{H^s}\|\partial_j\mathcal{U}\|_{H^{s-1}}\|\mathcal{U}\|_{H^s}\lesssim \varepsilon c(\|\mathcal{U}\|_{H^s},\|\mathfrak{b}\|_{H^s})\|\mathcal{U}\|_{H^s}^2.
\end{align*}
In both these commutator estimates, we made use of the fact that the constant diagonal matrix component trivially cancels out in the commutator.

Finally, attention has to be paid to the source term too, since for instance now we can't apply the ODE (\ref{eq_systemorder2bo2}) to treat the original right hand side of the first equation due to the presence of the operator $\Lambda^s$. Thus we are left with
\begin{equation}\label{estimate_firstrighthand}
\begin{aligned}
\int_{\mathbb{R}^d}\Lambda^s\zeta\Lambda^s\nabla_x\mathfrak{b}\left(x- X_S\right)\,dx\cdot\dot{X}_S&\leqslant\sqrt{\varepsilon}\frac{|\dot{X}_S|}{\sqrt{\varepsilon}}\|\zeta\|_{H^s}\|\nabla_x\mathfrak{b}\|_{H^s}\\
&\leqslant \sqrt{\varepsilon}c(\|\mathcal{U}\|_{H^s})\left(\frac{1}{\varepsilon}|\dot{X}_S|^2+\|\nabla_x\mathfrak{b}\|_{H^s}^2\right).
\end{aligned}
\end{equation}
Here we remark that $|\dot{X}_S|^2\leqslant 2\varepsilon E_B^s$. Notice that it is at this point that we can no longer use the cancellation, thus loosing a smallness factor.

Moreover, we are required to estimate terms which involve the operator $\Lambda^s$ applied to a product, this is handled by the commutator estimates, giving us
\begin{equation*}
\varepsilon\int_{\mathbb{R}^d}\Lambda^s(\nabla_x\mathfrak{b}\left(x- X_S\right)\cdot\overline{V})\Lambda^s\zeta\,dx\leqslant\varepsilon c(\|\mathfrak{b}\|_{H^{s+1}})\|\mathcal{U}\|_{H^s}^2,
\end{equation*}
and with a simple upper bound, we have
\begin{equation*}
\frac{\mu}{2}\int_{\mathbb{R}^d}h\Lambda^s\nabla_x\left(\nabla_x^2\mathfrak{b}\left(x- X_S\right)\dot{X}_S\cdot\dot{X}_S\right)\cdot\Lambda^s\overline{V}\,dx\leqslant\mu c(\|\mathcal{U}\|_{H^s},\|\mathfrak{b}\|_{H^{s+3}})\left|\dot{X}_S\right|^2.
\end{equation*}

At last, just as with Lemma \ref{lem_derniersourceterm}. we can deduce the following estimate
\begin{equation}\label{estimate_justlikelemma}
\begin{aligned}
\mu \int_{\mathbb{R}^d}\int_0^t&h\Lambda^s\nabla_x^2\mathfrak{b}\left(x- X_S(\tau)\right)\ddot{X}_S(\tau)\cdot\Lambda^s\overline{V}\,d\tau\,dx\leqslant\mu c(\brck{\mathcal{U}}_{T,H^s},\|\mathfrak{b}\|_{H^{s+2}})(E_B^s(t)+E_B^s(0))\\
&+\sqrt{\varepsilon} c(c_{fric},\tilde{M}^{-1},\brck{\mathcal{U}}_{T,H^s},\|\mathfrak{b}\|_{H^{s+3}})\int_0^tE_B^s(\tau)\,d\tau+\mu c(\brck{\mathcal{U}}_{T,H^s},\|\mathfrak{b}\|_{W^{1,\infty}})t \|\nabla_x\mathfrak{b}\|_{H^{s+2}}^2.
\end{aligned}
\end{equation}

To sum it up, after an integration with respect to the time variable, with the definition of the energy functional $E_B^s(t)$ and the velocity estimate obtained from the $L^2$ estimate (Corollary \ref{cor_velocityestimate}.), we may write that
\begin{equation}
E_B^s(t)\leqslant 2E_B^s(0)+\sqrt{\varepsilon}\tilde{c}_{s} t\|\mathfrak{b}\|_{H^{s+3}}^2+\sqrt{\varepsilon}\tilde{c}_{s}\int_0^tE_B^s(\tau)\,d\tau,
\end{equation}
with the constant
\begin{equation*}
\tilde{c}_{s}= c(c_{fric},\tilde{M}^{-1},\brck{\mathcal{U}}_{T,H^s},\|\mathfrak{b}\|_{H^{s+3}}).
\end{equation*}

So we have the right terms in order to complete the estimate, again with Gr\"onwall's lemma.\quod

\subsection{Local in time existence theorem}
The energy estimate allows us to establish the main existence theorem for the coupled Boussinesq system, which states as follows
\begin{theo}\label{theo_longtimebous}
Let us consider the coupled system defined by equations (\ref{eq_systemorder2bous}). Let us suppose that for the initial value $\zeta_{in}$ and $\mathfrak{b}$ the lower bound condition (\ref{genminheightcond}) 
\begin{equation}
\exists h_{min}>0,\;\forall X\in\mathbb{R}^d,\;1+\varepsilon\zeta(X)-\varepsilon\mathfrak{b}(X)\geqslant h_{min}
\end{equation}
is satisfied. If the initial values $\zeta_{in}$ and $\overline{V}_{in}$ are in $\mathcal{X}^s(\mathbb{R}^d)$ with $s\in\mathbb{R}$, $s>d/2+1$, and $ V_{S_0}\in\mathbb{R}^d$ then there exists a maximal $T_0>0$ independent of $\varepsilon$ such that there is a unique solution 
\begin{equation*}
\begin{aligned}
&(\zeta,\overline{V})\in C\left(\left[0,\frac{T_0}{\sqrt{\varepsilon}}\right];\mathcal{X}^s(\mathbb{R}^d)\right)\cap C^1\left(\left[0,\frac{T_0}{\sqrt{\varepsilon}}\right];\mathcal{X}^{s-1}(\mathbb{R}^d)\right),\\ 
&X_S\in C^2\left(\left[0,\frac{T_0}{\sqrt{\varepsilon}}\right]\right)
\end{aligned}
\end{equation*}
with uniformly bounded norms for the system (\ref{eq_systemorder2bous}) with initial conditions $(\zeta_{in},\overline{V}_{in})$ and $(0,\sqrt{\varepsilon}V_{S_0})$.
\end{theo}

\prove For this demonstration we shall follow the footsteps of a classical Friedrichs type reasoning for (in general) symmetric hyperbolic systems, found for example in Chapter $16$ of \cite{taylorpde}. The reason for this has already been evoked in the previous section, an iterative scheme is not adapted to the nonlinear coupled Boussinesq system because it does not allow for the cancellation of the coupling terms in the energy estimates. With a carefully chosen Friedrichs smoothing of the equations, these cancellations can be preserved.

\noindent \textbf{1. A regularized system:} We shall first of all regularize the system with the help of the Friedrichs mollifier $J_\delta$.
\begin{defi}
For every $u\in L^2(\mathbb{R}^d)$ we have that for $\xi\in\mathbb{R}^d$
\begin{equation*}
\widehat{J_\delta u}(\xi)=\varphi(\delta\xi)\hat{u}(\xi),
\end{equation*}
with $\varphi$ a regular real valued even function defined on $\mathbb{R}^d$ with compact support, such that $\varphi(0)=1$.
\end{defi}
A slightly modified classical property of the mollifier entails the followings
\begin{lem}\label{lem_mollifierprop}
$(1)$ For every $s,t\in\mathbb{R}$, the operator $J_\delta$ acts from $\mathcal{X}^s$ onto $\mathcal{X}^t$, moreover there exists a constant $C(s,t,\delta)$ such that
\begin{equation*}
\|J_\delta u\|_{\mathcal{X}^t}\leqslant C(s,t,\delta)\|u\|_{\mathcal{X}^s}
\end{equation*}
for every $u\in\mathcal{X}^s$.
$(2)$ $J_\delta$ as a linear operator is continuous every $L^p(\mathbb{R}^d)$, $1\leqslant p\leqslant\infty$, furthermore for all $u\in L^p(\mathbb{R}^d)$
\begin{equation*}
\|J_\delta u\|_{L^p}\leqslant C\|u\|_{L^p}
\end{equation*}
with a constant $C$ independent of $\delta$.
\end{lem}

Using the mollifier, we propose the following symmetric regularized system
\begin{subequations}\label{eq_systemsymmorder2reg}
\begin{align}[left = {\,\empheqlbrace}]
\begin{split}
&S(J_\delta\mathcal{U}^\delta,X_S^\delta)D_{\mu}\partial_t\mathcal{U}^\delta + \sum_{j=1}^{d}J_\delta S(J_\delta\mathcal{U}^\delta,X_S^\delta)A_j(J_\delta\mathcal{U}^\delta,X_S^\delta)J_\delta\partial_j\mathcal{U}^\delta= \\
&\qquad\qquad\qquad =S(J_\delta\mathcal{U}^\delta,X_S^\delta)J_\delta B(J_\delta\mathcal{U}^\delta,X_S^\delta), 
\end{split}\label{eq_systemsymmorder2reg1}\\
&\ddot{X}_S^\delta(t)=\mathcal{F}[J_\delta\mathcal{U}_0^\delta]\left(t,X_S^\delta,\dot{X}_S^\delta\right),\label{eq_systemsymmorder2reg2}\\
&\mathcal{U}^\delta(0,\cdot)=\mathcal{U}_{in},\left(X_S^\delta,\dot{X}_S^\delta\right)(0)=(0,v_{S_0}).\nonumber 
\end{align}
\end{subequations}

Based on Lemma \ref{lem_mollifierprop}, we may deduce that the regularized system (\ref{eq_systemsymmorder2reg}) is in fact an ODE (in the Fourier space) on any $\mathcal{X}^s$ Banach-space, the regularization guarantees that the nonlinear operator on the right hand side in its canonical form is regular thus uniformly Lipschitz and continuous in time. So by the Picard--Lindel\"of theorem we may deduce that there exists a solution $\mathcal{U}^\delta\in C([0,T_\delta];\mathcal{X}^s)$ and $X_S^\delta\in C^2([0,T_\delta])$.

\noindent \textbf{2. A priori estimate for the regularization:} Following the steps of the a priori estimates proved in the previous section, the estimate in Proposition \ref{prop_aprioriestimatebous2} holds for our regularized system as well, since by the careful choice of regularization in system (\ref{eq_systemsymmorder2reg}) the cancellations are preserved. So we have that
\begin{equation}
\|\mathcal{U}^\delta\|_{\mathcal{X}^s}^2+\frac{1}{\varepsilon}\left|\dot{X}_S^\delta(t)\right|^2\leqslant e^{\lambda t}\left(\|\mathcal{U}_{in}\|_{\mathcal{X}^s}^2+\left|v_{S_0}\right|^2+\sqrt{\varepsilon}c_st\|\mathfrak{b}\|_{H^{s+3}}^2\right),
\end{equation}
with $\lambda=\sqrt{\varepsilon}c_{s}$, $t\in[0,T_\delta]$.

\noindent \textbf{3. Uniformization of the time interval:} Here the hypothesis $s>d/2+1$ is important since we want to make use of the Sobolev embedding $H^s\hookrightarrow W^{1,\infty}$.
\begin{lem}
The regularized problem (\ref{eq_systemsymmorder2reg}) has a solution on $[0,\varepsilon^{-1/2}T_0]$ with $T_0$ independent of $\delta$ and $\varepsilon$.
\end{lem}
We have an estimate of the form
\begin{equation*}
\frac{d}{dt}E_B^s(\mathcal{U}^\delta,X_S^\delta)(t)\leqslant \sqrt{\varepsilon}c(c_{fric},\tilde{M}^{-1},\brck{\mathcal{U}^\delta}_{T,H^s},\|\mathfrak{b}\|_{H^{s+3}})\left(E_B^s(\mathcal{U}^\delta,X_S^\delta)(t)+\|\mathfrak{b}\|_{H^{s+3}}^2\right),
\end{equation*}
just before using Gr\"onwall's lemma in the higher order energy estimates. By a change of variable of the time parameter of the form $t=\varepsilon^{-1/2}t'$, it is clear to see that we have an uniform upper bound for a long time regime (with the variable $t'$). This implies that solutions to the regularized system (\ref{eq_systemsymmorder2reg}) exist for a time $\varepsilon^{-1/2}T_{0,\delta}$ with $T_{0,\delta}$ independent of $\varepsilon$.

Furthermore, we have that for every $\tilde{T}>0$, $t\in\left[0,\min\{\tilde{T},\varepsilon^{-1/2}T_{0,\delta}\}\right]$
\begin{equation}
\frac{d}{dt}\left(\|\mathcal{U}^\delta(t,\cdot)\|_{\mathcal{X}^s}^2+\frac{1}{\varepsilon}\left|\dot{X}_S^\delta(t)\right|^2\right)\leqslant F\left(\|\mathcal{U}^\delta(t,\cdot)\|_{\mathcal{X}^s}^2+\frac{1}{\varepsilon}\left|\dot{X}_S^\delta(t)\right|^2\right)
\end{equation}
with $F$ being a regular function independent of $\delta$. By the Picard--Lindel\"of theorem, there exists $T_0>0$ $(\varepsilon^{-1/2}T_0\leqslant \tilde{T})$ such that the ordinary differential equation
\begin{align*}[left = {\,\empheqlbrace}]
&y'(t)=F(y(t))\\
&y(0)=\|\mathcal{U}_{in}\|_{\mathcal{X}^s}^2+\left|V_{S_0}\right|^2
\end{align*}
has a unique solution on $[0, \varepsilon^{-1/2}T_0]$. By Gr\"onwall's lemma and a standard comparison theorem for ODEs we deduce that
\begin{equation*}
\|\mathcal{U}^\delta(t,\cdot)\|_{\mathcal{X}^s}^2+\frac{1}{\varepsilon}\left|\dot{X}_S^\delta(t)\right|^2\leqslant y(t).
\end{equation*}
\quod

\noindent \textbf{4. Convergence:} Let us define the following supplementary function space
\begin{equation*}
E_{T_0}^s=L^\infty\left(\left[0,\frac{T_0}{\sqrt{\varepsilon}}\right];\mathcal{X}^s\right)\cap W^{1,\infty}\left(\left[0,\frac{T_0}{\sqrt{\varepsilon}}\right];\mathcal{X}^{s-1}\right).
\end{equation*}

Solutions $\mathcal{U}^\delta$ of the regularized problem clearly belong to $E_{T_0}^s$ with $X_S^\delta\in W^{1,\infty}$. Thus, the family $\{\mathcal{U}^\delta\}_\delta$ is bounded in $E_{T_0}^s$ so it has a weakly convergent subsequence in $E_{T_0}^s$ towards a function $\mathcal{U}\in E_{T_0}^s$. 

Since the inclusion $H_{loc}^s(\mathbb{R}^d)\hookrightarrow H_{loc}^{s-1}(\mathbb{R}^d)$ is compact, by the Arzela--Ascoli theorem we may extract a strongly convergent subsequence from it in $C([0, \varepsilon^{-1/2}T_0];\mathcal{X}_{loc}^{s-1})$ locally.

By interpolation inequalities we have that $\mathcal{U}\in C^\sigma([0, \varepsilon^{-1/2}T_0];\mathcal{X}_{loc}^{s-\sigma})$ for each $\sigma\in(0,1)$. Moreover, since the inclusion $H^{s-\sigma}(\mathbb{R}^d)\hookrightarrow C^1(\mathbb{R}^d)$ is also compact for sufficiently small $\sigma>0$, we may deduce that the subsequence is converging in $C([0, \varepsilon^{-1/2}T_0];C_{loc}^1(\mathbb{R}^d))$, with $X_S^\delta$ converging in $C^1[0,\varepsilon^{-1/2}T_0]$. With this subsequence we shall no problem in passing to the limit in $\delta$ for the regularized system (\ref{eq_systemsymmorder2reg}) in the sense of distributions, leading to a solution of the problem in $E_{T_0}^s$.

\noindent \textbf{5. Additional regularity:} In fact, by being more careful with the estimates, we may deduce that the solution $\mathcal{U}$ is in
\begin{equation*}
C\left(\left[0,\frac{T_0}{\sqrt{\varepsilon}}\right];\mathcal{X}^s\right)\cap C^1\left(\left[0,\frac{T_0}{\sqrt{\varepsilon}}\right];\mathcal{X}^{s-1}\right).
\end{equation*}
Essentially, the main idea is to prove that the norm $\|\mathcal{U}(t,\cdot)\|_{\mathcal{X}^s}$ is in fact a (Lipschitz-)continuous function of $t$, since it is the limit of $\|J_\delta\mathcal{U}(t,\cdot)\|_{\mathcal{X}^s}$. For more details, we refer to Chapter 16 of \cite{taylorpde}.

\noindent \textbf{6. Uniqueness:} By taking the difference of two solutions for the system (\ref{eq_systemorder2bous}), they consequently satisfy a similar system, thus by the previous a priori estimate, with $0$ right hand side, we may conclude that this difference has to be $0$ as well. \quod 

This concludes the proof of the well-posedness theorem concerning the coupled fluid-solid system in the Boussinesq regime (\ref{eq_systemorder2bous}). As one can clearly see from the demonstration, Remark \ref{rem_existencetimevelo}. on the nature of the time of existence stays valid, so solutions are guaranteed over a time of order $\mathcal{O}(\varepsilon^{-1/2})$.

\subsection{Towards a more refined solid model}\label{sec_lastsection}

As we remarked in the beginning of the analysis of the Boussinesq regime, the solid equation was consistent with the full Newton equation only at $\mathcal{O}(\sqrt{\mu})$ (Proposition \ref{prop_consistorder2}.). In order to be more consistent with the equation, we have to continue the asymptotic development of $\Phi$ with respect to $\mu$, since the loss of consistency is due to the integral term
\begin{equation*}
-\frac{c_{fric}}{\tilde{M}\mu^{3/2}}\frac{\dot{X}_S}{|\dot{X}_S|+\overline{\delta}}\int_{I(t)}P_{z=-1+\mu b}\,dx.
\end{equation*}

For this, let us briefly elaborate an additional term in the asymptotic development. By Proposition 3.37. of \cite{lannesbible} we have that
\begin{equation}
\Phi(x,z) = \psi(x) -\mu\left(\frac{z^2}{2}+z\right)\Delta_x\psi + \mu z\partial_tb +\mathcal{O}(\mu^2),
\end{equation}
which gives us
\begin{equation*}
\varepsilon\partial_t\Phi = \varepsilon\partial_t\psi(x) -\varepsilon\mu\left(\frac{z^2}{2}+z\right)\Delta_x\partial_t\psi + \varepsilon\mu z\partial_{tt}b +\mathcal{O}(\varepsilon\mu^2),
\end{equation*}
as well as
\begin{equation*}
\frac{\varepsilon^2}{2}|\nabla_x\Phi|^2 = \frac{\varepsilon^2}{2}|\nabla_x\psi|^2+\mathcal{O}(\varepsilon^2\mu),\quad\frac{\varepsilon^2}{2\mu}|\partial_z\Phi|^2 =\mathcal{O}(\varepsilon^2\mu).
\end{equation*}
This means that the pressure formula (\ref{form_pressureassymptotic}) takes the form
\begin{equation*}
P(x,z) = \frac{P_{atm}}{\varrho gH_0} - z +\varepsilon\zeta - \varepsilon\mu\left(\frac{z^2}{2}+z\right)\Delta_x\zeta - \varepsilon\mu z\partial_{tt}b + \mathcal{O}(\varepsilon\mu^2),
\end{equation*}
so an evaluation at the bottom gives
\begin{equation}\label{form_pressionameliore}
P_{bott} = \frac{P_{atm}}{\varrho gH_0} + h + \frac{1}{2}\varepsilon\mu\Delta_x\zeta + \varepsilon\mu\partial_{tt}b + \mathcal{O}(\varepsilon\mu^2),
\end{equation}
and as such its integral over the support of the bottom ($I(t)$) is
\begin{equation*}
\int_{I(t)}P_{z=-1+\mu b}\,dx = \tilde{M}(\tilde{c}_{solid} - \varepsilon)+\varepsilon\int_{I(t)}\zeta\,dx +\frac{\varepsilon\mu}{2}\int_{I(t)}\Delta_x\zeta\,dx + \mathcal{O}(\varepsilon\mu^2).
\end{equation*}

By keeping the approximation of the pressure for the pressure term in the Newton's equation (\ref{eq_ndnewton}), we can recover a solid equation consistent at order $\mathcal{O}(\mu^{3/2})$ of the form
\begin{equation}\label{eq_newtonorder3bous}
\ddot{X}_S=-\frac{c_{fric}}{\sqrt{\mu}}\left(\frac{1}{\varepsilon}\tilde{c}_{solid}+\frac{1}{\tilde{M}}\int\limits_{\operatorname{supp}(\mathfrak{b})+ X_S}\zeta\,dx+\frac{\mu}{2\tilde{M}}\int\limits_{\operatorname{supp}(\mathfrak{b})+ X_S}\Delta_x\zeta\,dx\right)\frac{\dot{X}_S}{\left|\dot{X}_S\right|+\overline{\delta}}+\frac{\varepsilon}{\tilde{M}}\int_{\mathbb{R}^d}\zeta\nabla_x\mathfrak{b}(x-X_S)\,dx.
\end{equation}

One may obtain the same results for this model as the ones presented in Theorem \ref{theo_longtimebous}.

We would like to point out one particularity of the aforementioned computations. The integral of $\varepsilon\mu\partial_t^2b$ disappeared due to the fact that $b$ is smooth and of support compact. If one were to use the refined pressure formula (\ref{form_pressionameliore}) to compute the integral in the pressure term of Newton's equation as well, one would find an additional nonzero term, namely
\begin{equation}
\varepsilon\mu\int_{\mathbb{R}^d}(\nabla_xb\cdot\ddot{X}_S)\nabla_xb\,dx =: \mathcal{M}_{\nabla_xb}\ddot{X}_S,
\end{equation}
where the linear map $\mathcal{M}_{\nabla_xb}$ can be represented as a matrix, that is, in addition, positive semi-definite. In fact $\mathcal{M}_{\nabla_xb}$ stands for the so called added mass effect, or virtual mass effect, corresponding to an added inertia due to the solid accelerating/decelerating in the fluid medium, thus deflecting/moving some volume of the surrounding fluid as well. As one can see, in the weakly nonlinear Boussinesq regime (\ref{hyp_boussinesqregime}) this term was not present since it was of order $\mathcal{O}(\mu^2)$, however if one were to study the general second order asymptotic regime (meaning the Serre--Green--Naghdi equations), that is without the additional assumption of $\varepsilon=\mathcal{O}(\mu)$, this term (and many other nonlinearities) would be present.

\section*{Conclusion}

In the present paper, we established a coupled physical model of the water waves problem with a freely moving object on the bottom of the fluid domain. We deduced the exact coupled system and analyzed two different shallow water asymptotic regimes (with respect to the shallowness parameter $\mu$): the nonlinear Saint-Venant system and the Boussinesq system. We established local in time existence results as well as a uniqueness theorem for both cases and we improved the existence time for the weakly nonlinear Boussinesq regime.

Another possible approach would be to consider the full Green--Naghdi system for the $\mathcal{O}(\mu^2)$ asymptotic regime and establish the coupled system and possibly well-posedness results for it. This would yield a non-hydrostatic pressure formula, and consequently a more complex equation for the solid motion, incorporating the added mass effect, briefly elaborated in the last section.

An even more general scenario can be envisioned, that is to handle the full problem formulated in the first section, treating the coupled problem (\ref{eq_ndzakharovmovingbottom}) with (\ref{eq_ndnewton}).

To complement the theoretical results, a numerical study is to follow this article in order to verify the applicability of the system as well as to compare it with other existing methods to treat wave-structure interaction problems (\cite{zoazoazoo}, \cite{abadibadou}).

\bibliographystyle{alpha}
\bibliography{biblio}

\newcommand{\etalchar}[1]{$^{#1}$}
\begin{thebibliography}{ACDNn17}

\bibitem[ABZ11]{alazardbottommove}
Thomas Alazard, Nicolas Burq, and Claude Zuily.
\newblock On the {C}auchy problem for the water waves with surface tension.
\newblock {\em Duke Mathematical Journal}, 158:413--499, 2011.

\bibitem[ABZ14]{alazardgravitycauchy}
Thomas Alazard, Nicolas Burq, and Claude Zuily.
\newblock On the {C}auchy problem for gravity water waves.
\newblock {\em Inventiones Mathematicae}, 198:71--163, 2014.

\bibitem[ACDNn17]{abadibadou}
St\'ephane Abadie, Marcela Cruchaga, Benoit Ducassou, and Jonathan Nu\~nez.
\newblock A fictitious domain approach based on a viscosity penalty method to
  simulate wave/structure interaction.
\newblock {\em Journal of Hydraulic Research}, pages 1--16, 2017.

\bibitem[AELS14]{LehmannTheWC}
Reza Alam, Ryan Elandt, Marcus Lehmann, and Mostafa Shakeri.
\newblock {T}he {W}ave {C}arpet: Development of a {S}ubmerged {P}ressure
  {D}ifferential {W}ave {E}nergy {C}onverter.
\newblock 2014.
\newblock 30th Symposium on Naval Hydrodynamics, Hobart, Australia.

\bibitem[AG07]{alinhacgerard}
Serge Alinhac and Patrick G{\'e}rard.
\newblock {\em Pseudo-differential operators and the {N}ash--{M}oser theorem},
  volume~82 of {\em Graduate Studies in Mathematics}.
\newblock American Mathematical Society, 2007.

\bibitem[ASL08a]{alvarezlanneslong}
Borys Alvarez-Samaniego and David Lannes.
\newblock Large time existence for 3d water-waves and asymptotics.
\newblock {\em Inventiones Mathematicae}, 171:485--541, 2008.

\bibitem[ASL08b]{alvarezlannesGN}
Borys Alvarez-Samaniego and David Lannes.
\newblock A {N}ash--{M}oser theorem for singular evolution equations.
  {A}pplications to the {S}erre and {G}reen--{N}aghdi equations.
\newblock {\em Indiana University Mathemaics Journal}, 57:97--131, 2008.

\bibitem[BCL05]{davidboussinesq}
Jerry~Lloyd Bona, Thierry Colin, and David Lannes.
\newblock Long wave approximation for water waves.
\newblock {\em Archive for Rational Mechanics and Analysis}, 178:373--410,
  2005.

\bibitem[BCS04]{bonachen}
Jerry~Lloyd Bona, Min Chen, and Jean-Claud Saut.
\newblock {B}oussinesq equations and other systems for small-amplitude long
  waves in nonlinear dispersive media. {II}: {T}he nonlinear theory.
\newblock {\em Nonlinearity}, 17:925--952, 2004.

\bibitem[Ber06]{berthelotmechanic}
Jean-Marie Berthelot.
\newblock {\em M\'ecanique des solides rigides}.
\newblock Tec and Doc Lavoisier, 2006.

\bibitem[Bur16]{cosminboussinesq}
Cosmin Burtea.
\newblock New long time existence results for a class of {B}oussinesq-type
  systems.
\newblock {\em Journal de Math\'ematiques Pures et Appliqu\'ees},
  109(2):203--236, 2016.

\bibitem[Cha07]{chazelbottominfluence}
Florent Chazel.
\newblock Influence of bottom topography on long water waves.
\newblock {\em ESAIM: Mathematical Modeling and Numerical Analysis},
  41(4):771--799, 2007.

\bibitem[Che03]{chengeneralboussinesq}
Min Chen.
\newblock Equations for bi-directional waves over an uneven bottom.
\newblock {\em Mathematics and Computers in Simulation}, 62:3--9, 2003.

\bibitem[CLS12]{craiglannesrough}
Walter Craig, David Lannes, and Catherine Sulem.
\newblock Water waves over a rough bottom in the shallow water regime.
\newblock {\em Annales de l'Institut Henri Poincar\'e (C) Non Linear Analysis},
  29:233--259, 2012.

\bibitem[CM06]{cottetmaitrenumerics}
Georges-Henri Cottet and Emmanuel Maitre.
\newblock A level set method for fluid-structure interactions with immersed
  surfaces.
\newblock {\em Mathematical Models and Methods in Applied Sciences},
  16:415--438, 2006.

\bibitem[Cra85]{craiglangrangian}
Walter Craig.
\newblock An existence theory for water waves and the {B}oussinesq and
  {K}orteweg--de {V}ries scaling limits.
\newblock {\em Communications in Partial Differential Equations},
  10(8):787--1003, 1985.

\bibitem[CS93]{craigsulem1}
Walter Craig and Catherine Sulem.
\newblock Numerical simulation of gravity waves.
\newblock {\em Journal of Computational Physics}, 108:73--83, 1993.

\bibitem[CSS92]{craigsulem2}
Walter Craig, Catherine Sulem, and Pierre-Louis Sulem.
\newblock Nonlinear modulation of gravity waves: a rigorous approach.
\newblock {\em Nonlinearity}, 5:497--522, 1992.

\bibitem[DNZ15]{zoazoazoo}
Denys Dutykh, Hayk Nersisyan, and Enrique Zuazua.
\newblock Generation of two-dimensional water waves by moving bottom
  disturbances.
\newblock {\em IMA Journal of Applied Mathematics}, 80(4):1235--1253, 2015.

\bibitem[dP16]{poyferreshore}
Thibault de~Poyferr\'e.
\newblock A priori estimates for water waves with emerging bottom.
\newblock 2016.
\newblock arXiv:1612.04103.

\bibitem[GIL{\etalchar{+}}14]{waveenergyconverters}
Mateus das~Neves Gomes, Li\'ercio~Andr\'e Isoldi, Max Letzow, Luiz
  Alberto~Oliveira Rocha, Elizaldo Domingues~dos Santos, Fl\'avio~Medeiros
  Seibt, and Jeferson~Avila Souza.
\newblock Computational modeling applied to the study of wave energy converters
  ({WEC}).
\newblock {\em Marine Systems and Ocean Technology}, 9:77--84, 2014.

\bibitem[GKSW95]{weikirbyboussinesq}
Stephan~T. Grilli, James~T. Kirby, Ravishankar Subramanya, and Ge~Wei.
\newblock A fully nonlinear {B}oussinesq model for surface waves. {P}art 1.
  {H}ighly nonlinear unsteady waves.
\newblock {\em Journal of Fluid Mechanics}, 294:71--92, 1995.

\bibitem[GN07]{guyennespectrum}
Philippe Guyenne and David~P. Nicholls.
\newblock A high-order spectral method for nonlinear water waves over moving
  bottom topography.
\newblock {\em SIAM Journal of Scientific Computing}, 30(1):81--101, 2007.

\bibitem[HI15]{iguchiGNtsunami}
Fujiwara Hiroyasu and Tatsuo Iguchi.
\newblock A shallow water approximation for water waves over a moving bottom.
\newblock {\em Advanced studies in pure mathematics}, 64:77--88, 2015.

\bibitem[Igu09]{iguchishallow}
Tatsuo Iguchi.
\newblock A shallow water approximation for water waves.
\newblock {\em Journal of Mathematics of Kyoto University}, 49(1):13--55, 2009.

\bibitem[Igu11]{iguchibottommove}
Tatsuo Iguchi.
\newblock A mathematical analysis of tsunami generation in shallow water due to
  seabed deformation.
\newblock {\em Proceedings of the Royal Society of Edinburgh Section A:
  Mathematics}, 141:551--608, 2011.

\bibitem[KN86]{kanonishida}
Tadayoshi Kano and Takaaki Nishida.
\newblock A mathematical justification for {K}orteweg--de {V}ries equation and
  {B}oussinesq equation of water surface waves.
\newblock {\em Osaka Journal of Mathematics}, 23(2):389--413, 1986.

\bibitem[Lan05]{lanneswp}
David Lannes.
\newblock Well-posedness of the water-waves equations.
\newblock {\em Journal of the American Mathematical Society}, 18(3):605--654,
  2005.

\bibitem[Lan13]{lannesbible}
David Lannes.
\newblock {\em The water waves problem: mathematical analysis and asymptotics},
  volume 188 of {\em Mathematical Surveys and Monographs}.
\newblock American Mathematical Society, 2013.

\bibitem[Lan17]{lannesbonneton}
David Lannes.
\newblock On the dynamics of floating structures.
\newblock {\em Annals of Partial Differential Equations}, 3(1), 2017.

\bibitem[Li06]{liapproximatif}
Yi~A. Li.
\newblock A shallow-water approximation to the full water wave problem.
\newblock {\em Communications on Pure and Applied Mathematics},
  59(9):1225--1285, 2006.

\bibitem[LM17]{lannesmetivier}
David Lannes and Guy M\'etivier.
\newblock The shoreline problem for the one-dimensional shallow water and
  {G}reen--{N}aghdi equations.
\newblock 2017.
\newblock HAL Archives Ouvertes, <hal-01614321>.

\bibitem[Mel15]{melinandtsunami}
Benjamin Melinand.
\newblock A mathematical study of meteo and landlside tsunamis: the {P}roudman
  resonance.
\newblock {\em Nonlinearity}, 28:4037--4080, 2015.

\bibitem[Mit09]{mitsotakistsunami}
Dimitrios~E. Mitsotakis.
\newblock Boussinesq systems in two space dimensions over a variable bottom for
  the generation and propagation of tsunami waves.
\newblock {\em Mathematics and Computers in Simulation}, 80(4):860--873, 2009.

\bibitem[Ovs74]{ovsjannikov}
Lev~V. Ovsjannikov.
\newblock To the shallow water theory foundation.
\newblock {\em Archives of Mechanics}, 26:407--422, 1974.

\bibitem[Per67]{peregrineoriginal}
D.~Howell Peregrine.
\newblock Long waves on a beach.
\newblock {\em Journal of Fluid Mechanics}, 27(4):815--827, 1967.

\bibitem[PR83]{papanicolaou}
George~C. Papanicolaou and Rodolfo~R. Rosales.
\newblock Gravity waves in a channel with a rough bottom.
\newblock {\em Studies in Applied Mathematics}, 68:89--102, 1983.

\bibitem[SW00]{schneiderlongwavelimit}
Guido Schneider and C.~Eugene Wayne.
\newblock The long-wave limit for the water wave problem {I}. {T}he case of
  zero surface tension.
\newblock {\em Communications on Pure and Applied Mathematics},
  53(12):1475--1535, 2000.

\bibitem[SX12]{sautxulong}
Jean-Claude Saut and Li~Xu.
\newblock The {C}auchy problem on large time for surface waves {B}oussinesq
  systems.
\newblock {\em Journal de Math\'ematiques Pures et Appliqu\'ees},
  97(6):635--662, 2012.

\bibitem[Tay97]{taylorpde}
Michael~E. Taylor.
\newblock {\em Partial Differential Equations: Nonlinear Equations}, volume 117
  of {\em Applied Mathematical Sciences}.
\newblock Springer, 1997.

\bibitem[TW92]{wutengwaterwaves}
Michelle~H. Teng and Theodore~Y. Wu.
\newblock Nonlinear water waves in channels of arbitrary shape.
\newblock {\em Journal of Fluid Mechanics}, 242:211--233, 1992.

\bibitem[Wu87]{wusolitontrain}
Theodore~Y. Wu.
\newblock Generation of upstream advancing solitons by moving disturbances.
\newblock {\em Journal of Fluid Mechanics}, 184:75--99, 1987.

\bibitem[Wu97]{wuwellposedness2D}
Sijue Wu.
\newblock Well-posedness in {S}obolev spaces of the full water wave problem in
  2-{D}.
\newblock {\em Inventiones Mathematicae}, 130(1):39--72, 1997.

\bibitem[Wu99]{wuwellposedness3D}
Sijue Wu.
\newblock Well-posedness in {S}obolev spaces of the full water wave problem in
  3-{D}.
\newblock {\em Journal of the American Mathematical Society}, 12(2):445--495,
  1999.

\bibitem[Zak68]{zakharovformula}
Vladimir~E. Zakharov.
\newblock Stability of periodic water waves of finite amplitude on the surface
  of a deep fluid.
\newblock {\em Journal of Applied Mechanics and Technical Physics}, 9:190--194,
  1968.

\end{thebibliography}

\end{document}